# WAVELET REGRESSION IN RANDOM DESIGN WITH HETEROSCEDASTIC DEPENDENT ERRORS

By Rafał Kulik[1] and Marc Raimondo

*University of Ottawa and University of Sydney*

We investigate function estimation in nonparametric regression models with random design and heteroscedastic correlated noise. Adaptive properties of warped wavelet nonlinear approximations are studied over a wide range of Besov scales, $f \in \mathcal{B}^s_{\pi,r}$, and for a variety of $L^p$ error measures. We consider error distributions with Long-Range-Dependence parameter $\alpha, 0 < \alpha \leq 1$; heteroscedasticity is modeled with a design dependent function $\sigma$. We prescribe a tuning paradigm, under which warped wavelet estimation achieves partial or full adaptivity results with the rates that are shown to be the minimax rates of convergence. For $p > 2$, it is seen that there are three rate phases, namely the dense, sparse and long range dependence phase, depending on the relative values of $s, p, \pi$ and $\alpha$. Furthermore, we show that long range dependence does not come into play for shape estimation $f - \int f$. The theory is illustrated with some numerical examples.

**1. Introduction.**

1.1. *Random design regression with LRD errors.* Consider the random design regression model

(1.1) $$Y_i = f(X_i) + \sigma(X_i)\varepsilon_i, \qquad i = 1, \ldots, n,$$

where $X_i$'s are independent identically distributed (i.i.d.) random variables with a compactly supported density $g$, $\sigma(\cdot)$ is a deterministic function and

Received December 2008.
[1]Supported by a grant from the Natural Sciences and Engineering Research Council of Canada.

Unfortunately, the second author, Marc Raimondo, passed away before this article was completed.

*AMS 2000 subject classifications.* Primary 62G05; secondary 62G08, 62G20.

*Key words and phrases.* Adaptive estimation, nonparametric regression, shape estimation, random design, long range dependence, wavelets, thresholding, maxiset, warped wavelets.

$(\varepsilon_i)_{i \geq 1}$ is a stationary Gaussian sequence that is independent of the $X_i$'s. The long range dependence (LRD) of the $\varepsilon_i$'s is described by a linear structure

$$\varepsilon_i = \sum_{m=0}^{\infty} a_m \eta_{i-m}, \qquad a_0 = 1, \tag{1.2}$$

where $(\eta_i)_{i \in \mathbb{Z}}$ is an i.i.d. Gaussian sequence and $\lim_{m \to \infty} a_m m^{(\alpha+1)/2} = 1$, for $\alpha \in (0,1)$,

$$\mathrm{Var}\left(\sum_{i=1}^{n} \varepsilon_i\right) \sim c_\alpha n^{2-\alpha}, \tag{1.3}$$

where $c_\alpha$ is a finite and positive constant. In view of (1.3), the limit case $\alpha = 1$ can be thought of as similar to weakly dependent errors case.

1.2. *Prologue: linear regression.* Consider the regression model (1.1) with $f(x) = a + bx$ and correlated errors (1.2). We refer to Chapter 9 of [1], where the case $\sigma(x) \equiv 1$ is treated. The least squares (LS) estimator of $b$ is

$$\hat{b} = \frac{\sum_{i=1}^{n} X_i Y_i}{\sum_{i=1}^{n} X_i^2},$$

the asymptotic properties of $\hat{b} - b$ depend on $\sum_{i=1}^{n} X_i \sigma(X_i) \varepsilon_i$. Then,

$$\mathrm{Var}\left(\sum_{i=1}^{n} X_i \sigma(X_i) \varepsilon_i\right) = n \mathrm{E}(X_1^2 \sigma^2(X_1)) \mathrm{E}(\varepsilon_1^2)$$

$$+ (\mathrm{E}(\sigma(X_1) X_1))^2 \sum_{i \neq l}^{n} \mathrm{Cov}(\varepsilon_i, \varepsilon_l). \tag{1.4}$$

In this setting, the LS estimator is $\sqrt{n}$-consistent when the latter term is of order $O(n)$, which occurs if and only if

$$\mathrm{E}(\sigma(X_1) X_1) = 0. \tag{1.5}$$

For $\sigma(\cdot) \equiv 1$, this is always true when $\mathrm{E}(X_1) = 0$, and, if $\mathrm{E}(X_1) \neq 0$, it is enough to center shift the design variables $X_i - \bar{X}$, where $\bar{X} = \frac{1}{n} \sum_{i=1}^{n} X_i$. When $\sigma(\cdot) \not\equiv 1$, condition (1.5) is not necessarily fulfilled, even if $\mathrm{E}(X_1) = 0$. This is illustrated in the long range dependence literature. For example, [26] derived $\sqrt{n}$-consistency of a generalized least squares estimator when $\sigma(\cdot) \equiv 1$. Condition (1.5) appears in assumption 1 and Theorems 2.1 and 2.2 of [13]. This example suggests that, even in a simple parametric setting, statistical properties of LS estimators depend on the behaviour of $\sigma(\cdot)$ with respect to the design distribution. For example, if the design is uniformly distributed, $X_1 \sim \mathcal{U}[-1,1]$, then (1.5) is written as,

$$\int_{-1}^{1} \sigma(u) u \, du = 0, \tag{1.6}$$

which holds for any even function $\sigma$. Note, however, that, in practice, $\sigma$ is not observable.

1.3. *Background: nonparametric regression.* The model (1.2) with general error terms of the form $\sigma(X_t, \varepsilon_t)$ was considered in [5]. Asymptotic properties of the Nadaraya–Watson kernel estimator are found in [8], where $\varepsilon_i$'s are assumed to be a functional of LRD Gaussian random variables; in [9], with $\varepsilon_i$ as an infinite order moving average, and in [24], with the $X_i$'s possibly LRD, not necessarily independent of the $\varepsilon_i$'s. Local linear estimation using kernel method was studied in [22] and [23] and in case of FARIMA–GARCH errors in [1] and [2]. The corresponding results for density estimation were obtained in [4, 6, 15] and [28].

A general message from these papers is that the limiting behaviour of nonparametric estimators depends on a delicate balance between the smoothing parameter (e.g., bandwidth) and the long memory parameter $\alpha$. To be more specific, we quote the following result from [29], derived in (1.1) with $\sigma \equiv 1$:

$$(1.7) \qquad R_{n,2,g}(\mathcal{B}_{\pi,r}^s) \asymp n^{-\min(2s/(2s+1),\alpha)},$$

where $R_{n,p,g}(\mathcal{B}_{\pi,r}^s)$ denotes the minimax weighted $L^p$-risk over a Besov space $\mathcal{B}_{\pi,r}^s$,

$$(1.8) \qquad R_{n,p,g}(\mathcal{B}_{\pi,r}^s) := \inf_{\hat{f}} \sup_{f \in \mathcal{B}_{\pi,r}^s} \mathrm{E}\|f - \hat{f}\|_{L^p(g)}^p$$

with

$$\|f - h\|_{L^p(g)} = \left(\int_0^1 |f(x) - h(x)|^p g(x)\,dx\right)^{1/p}.$$

We refer to Section 2.2 for the precise definition of Besov spaces in terms of wavelet coefficients. Here, $s$ is related to the smoothness of the target function $f$, whereas $\pi$ and $r$ are scale parameters. In (1.7) we see that there is an elbow in the rate of convergence and, hence, that the best possible rate depends on the relative value of $s$ and $\alpha$. For small values of $\alpha$, LRD has a detrimental effect on the rates of convergence, whereas, for larger values of $\alpha$, we obtain the same rate as if the errors were independent. This is of importance in the development of adaptive tuning procedures since, in practice, neither $s$ nor $\alpha$ is known (note, however, that $\alpha$ can be estimated). While, for $\alpha = 1$, different data-driven methods (e.g., cross-validation, plug-in) have been implemented for choosing the bandwidth (see, e.g., [30]), for $\alpha < 1$, the effect of LRD may influence such procedures. We refer to [6] and [15] for detailed studies in the density case. We are not aware about such considerations in the random design regression setting, however, similar phenomena are anticipated.

Indeed, not many adaptive methods for curve estimation in the presence of long memory in errors are available. To the best of our knowledge, [12] is one of the few papers in this direction, where an orthogonal series estimator with adaptive stopping rule is shown to achieve the minimax rate, similar to that of (1.7), in the model (1.2) with $\sigma(\cdot) \equiv 1$. In [12], it was also noticed that the rate of convergence for shape estimation $f^* = f - \int f$ does not involve $\alpha$ and is the same as if the errors were independent. This observation was later confirmed by the minimax results of [29].

1.4. *Rates of convergence of wavelet estimators.* In this paper, we study adaptive function estimation in the model (1.2), performances of estimators are given with respect to various $L^p$, $p \geq 2$, error measures. Introducing the maximal risk

$$(1.9) \qquad R_{n,p,g}(\hat{f}_n, \mathcal{B}_{\pi,r}^s) := \sup_{f \in \mathcal{B}_{\pi,r}^s} \mathrm{E}\|f - \hat{f}_n\|_{L^p(g)}^p,$$

we consider nonlinear warped wavelet estimators of the form

$$\hat{f}_n(x) = \sum_{(j,k) \in \Lambda_1} \hat{\beta}_{j,k} \mathbb{I}\{|\hat{\beta}_{j,k}| \geq \lambda\} \psi_{j,k}(G(x)),$$

where $G(x) = \int_{-\infty}^{x} g(u)\,du$ is the design distribution function and $(\psi_{j,k})$ is a wavelet family with enough regularity. We show the statistical parameters $\hat{\beta}_{j,k}$ and tuning parameters $\lambda, \Lambda_1$ may be constructed independently of $s$ and to achieve near optimal results. Moreover, the tuning parameter $\lambda$ can be chosen independently of $\alpha$ as long as, for all $j \geq 0$, $k$,

$$(1.10) \qquad \mathrm{E}(\psi_{j,k}(G(X_1))\sigma(X_1)) = 0.$$

Note that, for $\sigma(\cdot) \equiv 1$, the condition (1.10) is always satisfied, since wavelets are orthogonal to constants (Haar family included). We note the similarity between condition (1.5) in the parametric setting and condition (1.10) in the nonparametric scenario.

Introducing rate exponents

$$(1.11) \qquad \alpha_D := \frac{2s}{2s+1}, \qquad \alpha_S := \frac{2(s - (1/\pi - 1/p))}{2(s - 1/\pi) + 1},$$

we will show that

$$R_{n,p,g}(\hat{f}_n, \mathcal{B}_{\pi,r}^s) \leq C n^{-p/2\gamma} (\log n)^\kappa,$$

where

$$(1.12) \quad \gamma = \begin{cases} \alpha_D, & \text{if } \alpha > \alpha_D \text{ and } s > \dfrac{p - \pi}{2\pi} \text{ (dense phase)}, \\ \alpha_S, & \text{if } \alpha > \alpha_S \text{ and } \dfrac{1}{\pi} < s < \dfrac{p - \pi}{2\pi} \text{ (sparse phase)}, \\ \alpha, & \text{if } \alpha \leq \min(\alpha_S, \alpha_D) \text{ (LRD phase)}, \end{cases}$$

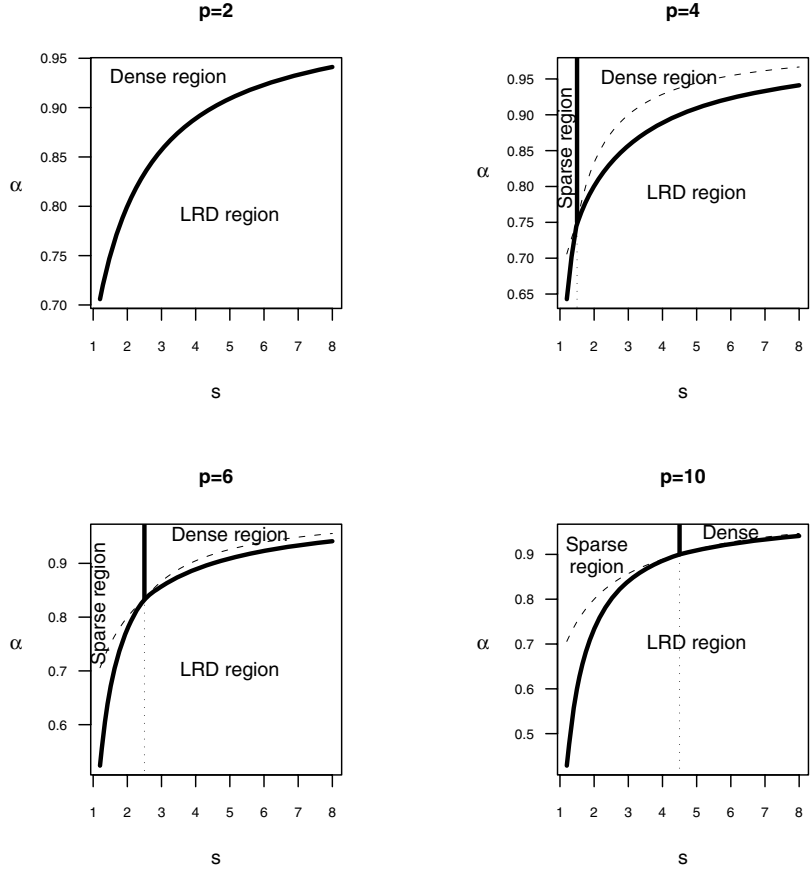

FIG. 1. *Illustration of the three rate phases (1.12) for Besov scales $\mathcal{B}_{1,1}^s$. Curve depicts rate exponent boundary, $\alpha(s) := \min(\frac{2s}{2s+1}, \frac{2(s-1+1/p)}{2s-1})$. Vertical line indicates the dense-s-parse boundary at $s = \frac{p-1}{2}$, for $p > 2$, this line meets the LRD boundary at $\alpha = \frac{p-1}{p}$.*

and $\kappa > 0$. This shows that convergence rates depends on the relative value of $\alpha$ with respect to $s$ but also on the relative value of $s$ with respect to $p$ and $\pi$.

We show, also, that our rates are optimal (up to a log term) in the minimax sense. Consequently, we generalize the result (1.7) to $p \geq 2$ and heteroscedastic errors. In particular, for $p = 2$ the rates agree with Yang's optimal rate, with a multiplicative log penalty, which is usual for adaptation. For $p > 2$ our results show that there are two elbows and three phases in the convergence rates, namely the dense phase, the sparse phase and the long range dependence phase. This is illustrated in Figure 1.

Furthermore, we show that, in the case of estimating the shape $f - \int f$, there is no LRD phase, which agrees with the previous findings in [12] and [29]. Finally, we will also show that, in the nonlinear wavelet estimator, we may replace $G(\cdot)$ with a corresponding empirical distribution function, and the resulting estimator still achieves the minimax rates.

## 2. Preliminaries.

2.1. *Warped wavelets.* Consider an orthonormal wavelet basis on the interval $I = [0, 1]$, $[\phi_{j,k}(x), \psi_{j,k}(x)]$, where $\phi$ denotes the scaling function and $\psi$ denotes the wavelet. Here, $j \geq 0, k = 0, \ldots, 2^j - 1$ and $\phi_{j,k}(x) = 2^{j/2}\psi(2^j x - k)$, $\psi_{j,k}(x) = 2^{j/2}\psi(2^j x - k)$. We refer to Chapter 7.5 of [21] for the construction of such a basis. For any function $f \in L^2[0,1]$, we have the following representation:

$$f(x) = \sum_{k=0}^{2^{j_0}-1} \alpha_{j_0,k}\phi_{j_0,k}(x) + \sum_{j=j_0}^{\infty}\sum_{k=0}^{2^j-1} \beta_{j,k}\psi_{j,k}(x), \tag{2.1}$$

where

$$\beta_{j,k} = \int_0^1 f(x)\psi_{j,k}(x)\,dx \tag{2.2}$$

denotes the wavelet coefficients associated to $f$, with the obvious corresponding definition for the scaling coefficients $\alpha_{j_0,k}$. The transformation (2.2) is called the wavelet transform (WT), and the representation (2.1) is called the inverse wavelet transform (IWT). In the case where $f$ is observed on a regular grid $i/n, i = 1, 2, \ldots, n$, both the WT and IWT can be computed in $O(n(\log n))$ steps using Mallat's pyramid algorithm. In the case where the function $f$ is observed along a random grid, the implementation of the standard WT (2.2) and IWT (2.1) requires some extra care.

A *warped* wavelet basis [19] is a *modified* wavelet basis representation specifically designed to handle random design regression model (1.2). The *modification* is suited to accommodate the design distribution function $G(\cdot)$. Provided that $f \circ G^{-1} \in L^2[0,1]$, we have the following representation:

$$f(x) = \sum_{k=0}^{2^{j_0}-1} \alpha_{j_0,k}\phi_{j_0,k}(G(x)) + \sum_{j=j_0}^{\infty}\sum_{k=0}^{2^j-1} \beta_{j,k}\psi_{j,k}(G(x)) \tag{2.3}$$

with

$$\beta_{j,k} = \int_0^1 f(x)g(x)\psi_{j,k}(G(x))\,dx, \tag{2.4}$$

where $g(x) = G'(x)$, and $\alpha_{j,k}$ is defined as in (2.4), with $\phi$ in place of $\psi$. By analogy with the standard case, we will refer to (2.4) and (2.3) as the warped wavelet transform (WWT) and the inverse warped wavelet transform (IWWT), respectively. Note that, by changing the variable in (2.4),

$$\beta_{j,k} = \int_0^1 f(G^{-1}(x))\psi_{j,k}(x)\,dx. \tag{2.5}$$

This shows that the WWT of $f$ is equivalent to the WT of $f \circ G^{-1}$. In the case where the function $f$ is observed along a random sequence $X_i$ with density $g$, the WWT and IWWT can be implemented in practice using a modification of Mallat's pyramid algorithm. This is further detailed in Section 5.1.

2.2. *Besov scales.* Throughout this paper, we will assume that $\psi$ is a compactly supported wavelet with $q, q > s$, vanishing moments and $\psi \in \mathbf{C}^q$ (see Chapter 9 of [21]). We further assume that the corresponding wavelet basis $(\psi_{j,k})$ satisfies the Temlyakov property as stated in [18]. Typical examples include the Daubechies wavelet family with $q$ vanishing moments. Finally, we consider wavelet basis on the interval $[0, 1]$ with appropriate boundary modifications (see [7]). In the light of (2.5), it is natural to express smoothness condition, with respect to $f \circ G^{-1}$ rather than $f$, as in [19]. We assume that $f \circ G^{-1} \in \mathcal{B}^s_{\pi,r}([0,1])$, where $s > \max\{\frac{1}{\pi}, \frac{1}{2}\}$. The latter condition may be written as $f \in L_\pi([0,1])$ and

$$f \circ G^{-1} = \sum_{j,k} \beta_{j,k}\psi_{j,k} \in \mathcal{B}^s_{\pi,r}(\mathcal{I})$$

$$\iff \sum_{j \geq 0} 2^{j(s+1/2-1/\pi)r}\left[\sum_{0 \leq k < 2^j} |\beta_{j,k}|^\pi\right]^{r/\pi} < \infty.$$

The parameter $s$ can be thought of as being related to the number of derivatives of $f$. With different values of $\pi$ and $r$, the Besov spaces capture a variety of smoothness features in a function, including spatially inhomogeneous behaviour.

**3. Minimax lower bounds over Besov balls.** In this section, we construct minimax lower bounds for the $L^p$ minimax risk, given in (1.8), for both dense and sparse case. As mentioned in Section 1.3, for the $L^2$-risk, homoscedastic errors and the dense case, the lower bound was obtained in [29].

To state our result, let us recall (1.11).

THEOREM 3.1. *Consider the model (1.2) and assume that $f \circ G^{-1} \in \mathcal{B}^s_{\pi,r}$. Furthermore, assume that $\inf_x \sigma(x) > 0$. Then, as $n \to \infty$,*

$$R_{n,p,g}(\mathcal{B}^s_{\pi,r}) \geq \begin{cases} C_p(n^{-p\alpha_D/2} \vee n^{-p\alpha/2}), & \text{if } s > \dfrac{p-\pi}{2\pi}; \\ C_p\left(\left(\dfrac{\log n}{n}\right)^{-p\alpha_S/2} \vee n^{-p\alpha/2}\right), & \text{if } \dfrac{1}{\pi} - \dfrac{1}{2} < s < \dfrac{p-\pi}{2\pi}, \end{cases}$$

where $C_p$ is a finite and positive constant. Furthermore, if $\Delta = \{f : \mathrm{E}[f(X_1)] = 0\}$, then

$$R_{n,p,g}(\mathcal{B}^s_{\pi,r} \cap \Delta) \geq \begin{cases} C_p(n^{-p\alpha_D/2}), & \text{if } s > \dfrac{p-\pi}{2\pi}; \\ C_p\left(\left(\dfrac{\log n}{n}\right)^{-p\alpha_S/2}\right), & \text{if } \dfrac{1}{\pi} - \dfrac{1}{2} < s < \dfrac{p-\pi}{2\pi}. \end{cases}$$

The above theorem means that if $f \in \Delta$, then the lower bounds are exactly the same as in the case of i.i.d. random errors. If this is not the case, then the rates are influenced by long memory. Furthermore, we can see that the dependence between the predictors and errors have no influence as long as $\sigma(\cdot)$ is bounded from below. Consequently, the theorem extends findings in [29] in several directions. First, it deals with $p \geq 2$; second, it identifies the elbow in the sparse case; third, it allows dependence between errors and predictors.

**4. Upper bounds for wavelets estimators.**

4.1. *Partial adaptivity.* By partial adaptivity, we mean that our estimator does not depend on $s$, but $G$ is known. Let

$$\Lambda_1 := \{(j,k), j_0 = -1 \leq j \leq j_1, k = 0, 1, \ldots, 2^j - 1\}$$

be the set of resolution levels. Here, the lowest resolution level $j_0 = -1$ corresponds to scaling contributions at resolution level $j = 0$ (i.e., $\psi_{-1,k} := \phi_{0,k}$ and $\beta_{-1,k} := \alpha_{0,k}$). The fine resolution level $j_1$ is set to be

(4.1) $$2^{j_1} \sim \frac{n}{\log n},$$

which is a classical condition. In practice, for a sample size $n$, the maximal number of resolution levels is set to be $2^{j_1} \sim n/2$; hence, condition (4.1) typically means that all resolution levels are used in (4.2).

The partially adaptive wavelet estimator we are going to consider is

(4.2) $$\hat{f}_n(x) = \sum_{(j,k) \in \Lambda_1} \hat{\beta}_{j,k} \mathbb{I}\{|\hat{\beta}_{j,k}| \geq \lambda\} \psi_{j,k}(G(x)),$$

(4.3) $$\hat{\beta}_{j,k} := \frac{1}{n} \sum_{i=1}^n \psi_{j,k}(G(X_i)) Y_i.$$

The theoretical level-dependent threshold parameter is set to be

(4.4) $$\lambda = \tau_0 \lambda_{n,j} := \tau_0(\tilde{\lambda}_n \vee \tilde{\lambda}_{n,j})$$
$$:= \tau_0 \left( \frac{\log n}{\sqrt{n}} \vee \mathbb{1}\{\mathrm{E}[\psi_{j,k}(G(X_1))\sigma(X_1)] \neq 0\} \frac{(\log n)^{1/2}}{n^{\alpha/2}} \right),$$

where $\tau_0$ is *large enough*. Note that, formally, the threshold depends on both $j$ and $k$; however, from theoretical point of view, $k$ is irrelevant. Furthermore, in simulation studies, we will average over all $k$ to get threshold depending on $j$ only.

The following theorem gives the convergence rates for the nonlinear wavelet estimator (4.2) according to the LRD index $\alpha$; recalling elbows location (1.12).

THEOREM 4.1. *Let $\hat{f}_n$ be the wavelet estimator (4.2) with (4.1), (4.3) and (4.4). Assume that $f \circ G^{-1} \in \mathcal{B}^s_{\pi,r}([0,1])$, $\pi \geq 1$, where $s > \max\{\frac{1}{\pi}, \frac{1}{2}\}$, and that $\sigma(\cdot)$ is bounded. Then,*

$$\mathrm{E}\|f - \hat{f}_n\|^p_{L^p(g)} \leq C n^{-p/2\gamma}(\log n)^\kappa,$$

*where*

$$\gamma = \begin{cases} \alpha_D, & \text{if } \alpha > \alpha_D \text{ and } s > \dfrac{p-\pi}{2\pi} \text{ (dense phase)}, \\ \alpha_S, & \text{if } \alpha > \alpha_S \text{ and } \dfrac{1}{\pi} < s < \dfrac{p-\pi}{2\pi} \text{ (sparse phase)}, \\ \alpha, & \text{if } \alpha \leq \min(\alpha_S, \alpha_D) \text{ (LRD phase)}, \end{cases}$$

*and $\kappa = p\gamma$ in the sparse and dense phase, $\kappa = 1$ in the LRD phase. If $\alpha = 1$, the LRD phase is not relevant.*

REMARK 4.2. When $\alpha = 1$, there is only one elbow on the convergence rate, provided that $p > 2$, switching from rate exponent $\alpha_D$ (dense phase) to rate exponent $\alpha_S$ (sparse phase). This is consistent with results obtained in the case of independent errors (see, e.g., [19]).

REMARK 4.3. For $\alpha < 1$ and $p > 2$, our rate results seems to be new, and we see that there is an additional elbow in the convergence rate switching from rate exponent $\alpha_D$ or $\alpha_S$ to $\alpha$ (LRD phase), depending on the relative value of $s$ and $\alpha$. This is illustrated in Figure 1. For $p = 2$, we note that there is only one elbow in the convergence rate, as we are either in the dense phase when $\alpha > \alpha_D$ or in the LRD phase when $\alpha \leq \alpha_D$. This is consistent with results of [12] and [29].

REMARK 4.4. Note that if, for all $j \geq 0$, $k = 0, \ldots, 2^j - 1$, the condition (1.10) holds, then the threshold (4.4) does not involve $\alpha$ (i.e., the estimator is constructed in the same way as if the errors were independent). The threshold (4.4) is then similar to the universal used in wavelet shrinkage (see, e.g., [10]). There is an additional multiplicative $(\log n)^{1/2}$ term, which is due to the martingale approximation of LRD sequences.

To gain some insight into condition (1.10), we note that, in the case of a uniform design distribution, this condition is written as, for all $j, k$

$$\int \psi_{j,k}(u)\sigma(u)\,du = 0,$$

which typically holds if $\sigma(u)$ is a polynomial function and $\psi$ has enough vanishing moments. Typically, this condition does not hold if $\sigma$ has some irregularities (jumps, cusps) or if $\sigma$ is oscillating at medium and high frequencies.

REMARK 4.5. Our definition (4.1) of $j_1$ is the same as the definition used in standard (nonwarped) estimation with independent errors. We note that it is less restrictive than the definition used in the warped wavelet estimation setting of [19]. Because of such choice of $j_1$, the bias is of smaller order than the bias in [19]. Consequently, in the sparse phase we have the restriction $s > 1/\pi$, as compared to $s > \frac{1}{2} + \frac{1}{\pi}$ in Proposition 2 of [19]. See also Remark 4.9.

REMARK 4.6. A comparison of our results with rate results obtained under a regular grid design, [17, 20] and [27], shows that randomization of the design improves rate performances. We illustrate this using the fixed design rate exponents, but similar inequalities hold in the sparse region. In the fixed design scenario, the dense region rate exponent is $\alpha s/(s + \alpha/2)$, which is always smaller than the exponent $\min\{(2s/(2s+1)), \alpha\}$ achievable under a random design.

REMARK 4.7. Using the weighted norm approximation of Theorem 4.1, we can conclude some results for the usual norm, even when $g(x_0) = 0$ for some $x_0 \in [0, 1]$. To see this, let $A = \{x \in [0, 1] : g(x) \neq 0\}$ and assume that the Lebesgue measure of $[0, 1] \setminus A$ is zero. If, now, $\|\cdot\|_p = \|\cdot\|_{L^p(1)}$ is the usual $L_p$-norm, then, with $1/q_1 + 1/q_2 = 1$, $q_1, q_2 > 1$, $l \in \mathbb{R}$,

$$\begin{aligned}
\mathrm{E}\|f - \hat{f}_n\|_p^p &= \int_A \mathrm{E}|f - \hat{f}_n|^p = \int_A \mathrm{E}|f - \hat{f}_n|^p \frac{g^l}{g^l} \\
&\leq \left(\int_0^1 \{\mathrm{E}|f - \hat{f}_n|^p\}^{q_1} g^{lq_1}\right)^{1/q_1} \left(\int_A g^{-lq_2}\right)^{1/q_2} \\
&\leq \left(\int_0^1 \mathrm{E}|f - \hat{f}_n|^{pq_1} g^{lq_1}\right)^{1/q_1} \left(\int_A g^{-lq_2}\right)^{1/q_2} \\
&= (\mathrm{E}\|f - \hat{f}_n\|_{L^{pq_1}(g)}^{pq_1})^{1/q_1} \left(\int_A g^{-lq_2}\right)^{1/q_2},
\end{aligned}$$

by choosing $lq_1 = 1$. Take, now, as in [19], $g(x) = (a + 1)x^a$, $x \in [0, 1]$. Then, the latter integral is finite as long as $a < (q_2 - 1)^{-1}$. On the other hand, we can apply Theorem 4.1 to conclude that in the dense and LRD phase the rates of convergence of $\mathrm{E}\|f - \hat{f}_n\|_p^p$ are the same as of $\mathrm{E}\|f - \hat{f}_n\|_{L^p(g)}^p$, as long as

$$(4.5) \qquad s > \frac{pq_1 - \pi}{2\pi} = \frac{p(1 + 1/(q_2 - 1) - \pi)}{2\pi}.$$

Note, however, that, if $a < (q_2 - 1)^{-1} < (2 + \pi - p)/p$, then

$$s > \frac{1}{\pi} > \frac{p(1 + 1/(q_2 - 1) - \pi)}{2\pi},$$

so that (4.5) becomes void. Consequently, for any $a < (2 + \pi - p)/p$, we can obtain the optimal rates. Of course, this approach does not work in the sparse case, because the resulting upper bound is not optimal (cf. Theorem 2 of [19]).

Furthermore, if $0 < m < g < M < \infty$, then the norms $\|\cdot\|_p$ and $\|\cdot\|_{L^p(g)}$ are equivalent.

4.2. *Full adaptivity.* By full adaptivity, we mean that our estimator does not depend on $s$ and $G$ is unknown. In this case, the fine resolution level $j_1$ in (4.1) has to be modified thusly:

$$(4.6) \qquad 2^{j_1} \sim \sqrt{\frac{n}{\log n}}.$$

In fact, in general (see Remark 7.6), we cannot use the same fine resolution level as in (4.1).

Assume that we have $2n$ observations from the model (1.2) coded as follows: the first $n$ observations are denoted by $X'_1, \ldots, X'_n$, the remaining as $X_1, \ldots, X_n$. The estimator that achieves the full adaptivity is

$$(4.7) \qquad \tilde{f}_n(x) = \sum_{(j,k) \in \Lambda_1} \tilde{\beta}_{j,k} \mathbb{I}\{|\tilde{\beta}_{j,k}| \geq \lambda\} \psi_{j,k}(\hat{G}_n(x)),$$

where, now, $\hat{G}_n$ is the empirical distribution function associated with $X'_1, \ldots, X'_n$ and

$$(4.8) \qquad \tilde{\beta}_{j,k} := \frac{1}{n} \sum_{i=1}^n \psi_{j,k}(\hat{G}_n(X_i)) Y_i.$$

THEOREM 4.8. *Consider the estimator (4.7) with (4.4), (4.6) and (4.8). Assume that $f \circ G^{-1} \in \mathcal{B}_{\pi,r}^s([0,1]) \cap \mathrm{Lip}_{1/2}$, $\pi \geq 1$, where $s > \max\{\frac{1}{\pi}, \frac{1}{2}\}$, and that $\sigma(\cdot)$ is bounded. Then, the rates of Theorem 4.1 remain valid with*

$$\gamma = \begin{cases} \alpha_D, & \text{if } \alpha > \alpha_D \text{ and } s > \dfrac{p - \pi}{2\pi} \text{ (dense phase)}, \\ \alpha_S, & \text{if } \alpha > \alpha_S \text{ and } \dfrac{1}{\pi} + \dfrac{1}{2} < s < \dfrac{p - \pi}{2\pi} \text{ (sparse phase)}, \\ \alpha, & \text{if } \alpha \leq \min(\alpha_S, \alpha_D) \text{ (LRD phase)}. \end{cases}$$

REMARK 4.9. Note that, in the sparse phase, there is the additional restriction $s > \frac{1}{\pi} + \frac{1}{2}$, as compared to Theorem 4.1. This is due to the larger bias, which, in turn, is due to choosing lower highest resolution level.

4.3. *Shape estimation.* As first noticed in [12], the effect of LRD is concentrated on the zero Fourier frequency component of the target function $f$ and corresponds to the scale $\int f$ of $f$. Keeping this in mind, it is possible to avoid (or reduce) the curse of LRD by considering the estimation of the shape of the function: $f - \int f$. Taking into account the design distribution in (2.4), we set

$$f^*(x) := f(x) - \int_0^1 f(G^{-1}(y))\, dy =: f(x) - c_{f,G}.$$

Note that the wavelet coefficient $\beta_{j,k}^*$ of $f^* \circ G^{-1}$ is equal to $\beta_{j,k}$. We set

$$(4.9) \qquad \hat{f}_n^* := \sum_{(j,k) \in \Lambda_1, j \neq -1} \hat{\beta}_{j,k} \mathbb{I}\{|\hat{\beta}_{j,k}| \geq \lambda\} \psi_{j,k}$$

and $\tilde{f}_n^*$, the corresponding fully adaptive estimator. The trick here is simply to remove the scaling coefficient. This is allowed, since $\int_0^1 f^*(G^{-1}(y))\, dy = 0$. In this way, there will be no LRD effect on the convergence rates.

THEOREM 4.10. *Let $\hat{f}_n^*$ be the wavelet estimator (4.9). Under assumptions of Theorem 4.1,*

$$(4.10) \qquad \mathrm{E}\|f^* - \hat{f}_n^*\|_p^p \leq C n^{-p/2\gamma} (\log n)^{p\gamma}.$$

*Under the assumptions of Theorem 4.8, the same bound is valid for $\tilde{f}_n^*$.*

Note that $\int_0^1 f^*(G^{-1}(y))\, dy = \mathrm{E}[f^*(X_1)] = 0$. Therefore, by comparing (4.10) with the second part of Theorem 3.1, we see that $f^*$ is estimated (up to a log term) with the optimal rates.

## 5. Finite sample properties.

5.1. *Implementation.* In our simulation studies, we focus on LRD effect. For this purpose, we assume that $U_{(1)} \leq U_{(2)} \leq \cdots \leq U_{(n)}$ denotes the ordered design sample from the uniform distribution, and $Y_{(1)}, \ldots, Y_{(n)}$ the corresponding observations of $Y_i$, not necessary ordered. If $\hat{G}_n$ is the empirical distribution function associated with $U_{(1)}, \ldots, U_{(n)}$, we have

$$(5.1) \qquad \frac{1}{n}\sum_{i=1}^n \psi_{j,k}(\hat{G}_n(U_{(i)})) Y_{(i)} = \frac{1}{n}\sum_{i=1}^n \psi_{j,k}(i/n) Y_{(i)}.$$

As noted in [3], in the case of a uniform design distribution, the ordered sample $U_{(1)}, \ldots, U_{(n)}$ may be used as a proxy for the regular grid $t_i = i/(n+1)$. Thus, in this case, (5.1) is computed by a simple application of Mallat's algorithm using the $Y_{(i)}$'s as input variables. This algorithm is implemented in the wavethresh R-package with various thresholding options, from which it is straightforward to compute function and shape estimators. This is the software (appropriately modified) we have used in the examples below.

*Data-based threshold.* As mentioned in Remark 4.4, if (1.10) holds, then the threshold is almost like in the usual fixed-design regression, with i.i.d. errors $\tau_0 \log n/\sqrt{n}$; here, with the additional log penalty. The parameter $\tau_0$ is estimated by a standard deviation of wavelet coefficients on the finest resolution level (option by.level=FALSE) or by computing standard deviation on each level separately (option by.level=TRUE).

The LRD part of the threshold may be chosen in the following way. First, note that $E[\psi_{j,k}(G(X_1))\sigma(X_1)]$ is just the wavelet coefficient of $\sigma(G^{-1}(\cdot))$. Therefore, we may perform a preliminary estimation and compute residuals, which serve as proxies for $\sigma(X_i)\varepsilon_i$. From this, we can estimate $\sigma(\cdot)$ and then the dependence index $\alpha$. If $\hat\sigma(\cdot)$ is the estimator of $\sigma(\cdot)$, then we may apply DWT to $\hat\sigma(\hat G_n^{-1}(i/n))$. Extracting the resulting wavelet coefficients on level $j$, we obtained the estimates of $E[\psi_{j,k}(G(X_1))\sigma(X_1)]$. For a given $j$, the level dependent threshold is obtained as the average over $k = 0, \ldots, 2^j - 1$.

5.2. *Examples.* We generate $Y_i$'s data according to (1.1) with Lidar, Bumps and Doppler target

$$(5.2) \qquad f(x) = (x(1-x))^{1/2} \sin\left(2\pi \frac{1.05}{x+1.05}\right),$$

a uniform design distribution $X_i = U_i \sim \mathcal{U}[0,1]$ and the following three $\sigma(\cdot)$ scenarios: (a) *homoscedastic scenario with* $\sigma(x) \equiv 0.1$ (constant noise level); (b) *heteroscedastic with* $\sigma(x) = 0.1\sqrt{\frac{12}{13}}(x+0.5)$ (linear noise level); and (c) *heteroscedastic with* $\sigma(x) = 0.1(\sin(\pi x) - \text{sign}(x-0.4))$ (irregular noise level). For calibration and comparison purposes, we quote, for scenario (a) with the Doppler target, the signal-to-noise ratio (SNR)

$$\text{SNR} = 10 \log_{10}\left(\frac{\int f^2}{\sigma^2}\right) \approx 9.34 \text{ (dB)}.$$

All other target function (Bumps and Lidar) were standardized to obtain the same SNR. Two different threshold parameters are considered, one given by (4.4) and the standard Donoho–Johnstone threshold. The noise level is estimated either on each level (option by.level=TRUE) or globally (option

TABLE 1
*Monte Carlo approximations to MSE of function estimator (4.2) of the `Doppler` target, with 1000 replications of the model (1.1), in scenario* (a), (b) *and* (c) *for some values of the dependence parameter d*

| | (a) | | (b) | | (c) | |
|---|---|---|---|---|---|---|
| $d$ | DJ thr | LRD thr | DJ thr | LRD thr | DJ thr | LRD thr |
| 0.000 | 0.0277 | 0.0277 | 0.0276 | 0.0305 | 0.0280 | 0.0329 |
| 0.150 | 0.0276 | 0.0276 | 0.02745 | 0.0288 | 0.0279 | 0.0319 |
| 0.300 | 0.0284 | 0.0284 | 0.0282 | 0.0287 | 0.0289 | 0.0315 |
| 0.325 | 0.0280 | 0.0280 | 0.0278 | 0.0281 | 0.0284 | 0.0316 |
| 0.350 | 0.0282 | 0.0282 | 0.0281 | 0.0282 | 0.0288 | 0.0319 |
| 0.375 | 0.0299 | 0.0299 | 0.0297 | 0.0299 | 0.0306 | 0.0335 |
| 0.400 | 0.0320 | 0.0320 | 0.0317 | 0.0319 | 0.0326 | 0.0350 |
| 0.425 | 0.0350 | 0.0350 | 0.0347 | 0.0347 | 0.0358 | 0.0383 |
| 0.450 | 0.0449 | 0.0449 | 0.0445 | 0.0446 | 0.0466 | 0.0486 |

by.level=FALSE). For such threshold values, we apply two threshold policies, Hard and Soft. Finally, Daubechies DB(6) and DB(2) wavelets are considered. For each of those scenarios we study the effect of the LRD parameter $\alpha$ on the performances of function estimator (4.2) and shape estimator (4.9) for sample sizes $n = 1024$.

Monte Carlo results for `Doppler` and `Bumps`, with $N = 1000$ replications and Daubechies DB(6) wavelet are summarized in Tables 1 and 2 on page 14. Notation DJ thr and LTD thr stands for Donoho–Johnstone universal threshold and the one given in (4.4), respectively.

The mean square error $\text{MSE} := \frac{1}{n} \sum_{i=1}^{n} (f(i/n) - \hat{f}_n(i/n))^2$ is plotted as a function of the dependence parameter as $d = (1 - \alpha)/2 \in (0, 1/2)$ in Figure 2. Here, $d$ corresponds to the fractional integration parameter as required to simulate LRD noise using `fracdiff` R-package.

*Analysis of the results.*

1. Figure 2 describes MSE, for the homoscedastic scenario (a). We can observe that the MSE seems to remain stable when the dependence is in the $[0, 0.35]$ range. Then, a sudden increase occurs after 0.35 suggesting that, for this simulated example, the LRD phase becomes active for very dependent error terms and confirming the detrimental effect of LRD in this region. This is also confirmed in Table 1. The similar effect is visible in the case of `Bumps` function, in Table 2.
2. We compare Donoho–Jonstone classical threshold with the one introduced in (4.4). Comparing left and right panels in Tables 1 and 2, we can see that there is completely no difference in case of the heteroscedastic noise. However, in case of the irregular noise [like scenario (c) above], the picture is not clear. In Doppler case the classical threshold performs better, on the other hand the level-dependent threshold (4.4) is preferable in case of Bumps target. This also applies to Lidar function.

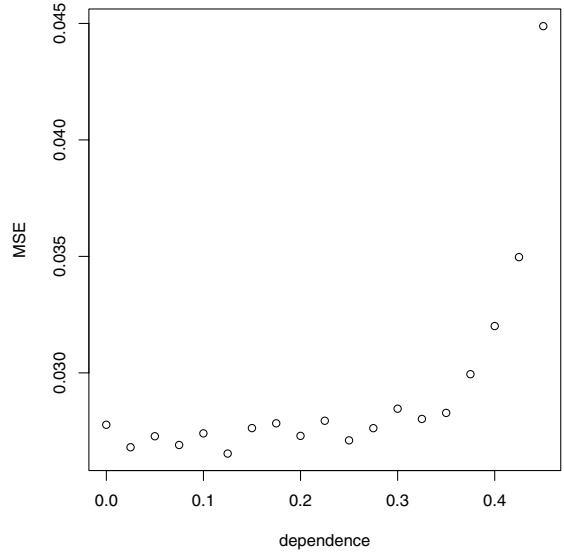

Fig. 2. *Monte Carlo approximation to MSE, $n = 1024$. Doppler function with $\sigma(x) = 0.1$.*

3. There is not too much difference between DB(2) and DB(6), as well as between Hard and Soft policy. However, the BY.LEVEL noise estimation

Table 2
*Monte Carlo approximations to MSE of function estimator (4.2) of the Bumps target, with 1000 replications of the model (1.1), in scenario* (a), (b) *and* (c) *for some values of the dependence parameter d*

| | (a) | | (b) | | (c) | |
|---|---|---|---|---|---|---|
| $d$ | DJ thr | LRD thr | DJ thr | LRD thr | DJ thr | LRD thr |
| 0.000 | 0.1295 | 0.1295 | 0.1293 | 0.1273 | 0.1297 | 0.1239 |
| 0.150 | 0.1298 | 0.1298 | 0.1297 | 0.1288 | 0.1301 | 0.1256 |
| 0.300 | 0.1297 | 0.1297 | 0.1294 | 0.1295 | 0.1300 | 0.1263 |
| 0.325 | 0.1301 | 0.1301 | 0.1297 | 0.1296 | 0.1308 | 0.1281 |
| 0.350 | 0.1309 | 0.1309 | 0.1306 | 0.1306 | 0.1315 | 0.1289 |
| 0.375 | 0.1328 | 0.1328 | 0.1324 | 0.1324 | 0.1334 | 0.1308 |
| 0.400 | 0.1340 | 0.1340 | 0.1335 | 0.1335 | 0.1349 | 0.1327 |
| 0.425 | 0.1377 | 0.1377 | 0.1372 | 0.1372 | 0.1389 | 0.1367 |
| 0.450 | 0.1462 | 0.1462 | 0.1456 | 0.1456 | 0.1487 | 0.1460 |

(i.e., estimation of $\tau_0 = \tau_{0,j}$) gives worse results in terms of MSE. The reason for this could be the following: variance estimator in case of LRD has slower rates of convergence then in the associated i.i.d. sequence. Consequently, on low frequencies (this is where LRD comes into play), the noise level estimates may not be very precise. The practical message is that, in LRD case, we should use the noise level estimates based on the highest resolution level.

**6. Proofs: lower bounds.** To obtain the lower bounds, we follow closely the ideas of [25]. Let us first introduce some notation. Denote $\mathbf{Y} = (Y_1, \ldots, Y_n)'$, $\varepsilon = (\varepsilon_1, \ldots, \varepsilon_n)'$, $\mathbf{1} = (1, \ldots, 1)'$, and, for any function $f$, let

$$f(\mathbf{X}) := (f(X_1), \ldots, f(X_n))$$

and $f(\mathbf{X})/\sigma(\mathbf{X})$ and $f(\mathbf{X}) * \sigma(\mathbf{X})$ be the coordinatewise division and multiplication, respectively, of two vectors. Furthermore, $\Xi$ is the covariance matrix of $\varepsilon$. With a slight abuse of notation, let $\Xi = (\xi_{il})_{i,l=1,\ldots,n}$ and $\Xi^{-1} = (\xi_{il}^{-1})_{i,l=1,\ldots,n}$ [of course, $(\xi_{il})^{-1} \neq \xi_{il}^{-1}$, in general].

For two functions $f, f_0$, denote by $\Lambda_n(f_0, f)$ the likelihood ratio

$$\Lambda_n(f_0, f) = d\mathbb{P}_{Y(f_0)}/d\mathbb{P}_{Y(f)},$$

where $\mathbb{P}_{Y(f)}$ is the distribution of the process $\{Y_i, i \geq 1\}$ when $f$ is true.

Note that the model (1.2) can be written as $\mathbf{Y}' = f(\mathbf{X})' + (\sigma(\mathbf{X}) * \varepsilon)'$. Then, we have, under $\mathbb{P}_{Y(f)}$,

$$
\begin{aligned}
2\ln \Lambda_n(f_0, f) \\
&= \left(\frac{\mathbf{Y} - f(\mathbf{X})}{\sigma(\mathbf{X})}\right)' \Xi^{-1} \left(\frac{\mathbf{Y} - f(\mathbf{X})}{\sigma(\mathbf{X})}\right) \\
&\quad - \left(\frac{\mathbf{Y} - f_0(\mathbf{X})}{\sigma(\mathbf{X})}\right)' \Xi^{-1} \left(\frac{\mathbf{Y} - f_0(\mathbf{X})}{\sigma(\mathbf{X})}\right) \\
&= -\left(\frac{f_0(\mathbf{X}) - f(\mathbf{X})}{\sigma(\mathbf{X})}\right)' \Xi^{-1} \left(\frac{(f_0(\mathbf{X}) - f(\mathbf{X}))}{\sigma(\mathbf{X})}\right) \\
&\quad + 2\left(\frac{f_0(\mathbf{X}) - f(\mathbf{X})}{\sigma(\mathbf{X})}\right)' \Xi^{-1} \varepsilon.
\end{aligned}
$$

(6.1)

In what follows, $\pi_0$ and $C_1, \ldots, C_4, C_p$ will be fixed and positive numbers.

*Sparse case.* This is the case when the hardest function to estimate is represented by one term in the wavelet expansion only. In this case, we use the result of Korostelev and Tsybakov (see [16], Lemma 10.1).

LEMMA 6.1. *Let $V$ be a functional space, and let $d(\cdot, \cdot)$ be a distance on $V$. Let $V$ contain the functions $f_0, f_1, \ldots, f_K$, such that:*

(a) $d(f_k, f_{k'}) \leq \delta > 0$ for $k = 0, 1, \ldots, K$, $k \neq k'$,
(b) $K \geq \exp(\lambda_n)$ for some $\lambda_n > 0$,
(c) $\ln \Lambda_n(f_0, f_k) = u_{nk} - v_{nk}$, where $v_{nk}$ are constants and $u_{nk}$ is a random variable such that for some $\pi_0 > 0$ we have $P_{f_k}(u_{nk} > 0) \geq \pi_0$,
(d) $\sup_k v_{nk} \leq \lambda_n$.

Then, for an arbitrary estimator $\hat{f}_n$,
$$\sup_{f \in V} \mathbb{P}_{Y_n^{(f)}}(d(\hat{f}_n, f) \geq \delta/2) \geq \pi_0/2.$$

To use this lemma, let us now choose $V = \{f_{jk} : 0 \leq k \leq 2^j - 1\}$, where $f_{jk}(x) = \beta_{j,k}\psi_{j,k}(G(x))$ [i.e., $f_{jk}(G^{-1}(u)) = \beta_{j,k}\psi_{j,k}(u)$, $f_0 \equiv 0$]. Since $f \circ G^{-1} \in \mathcal{B}_{\pi,r}^s$, we have $\beta_{j,k} \leq A 2^{-js'}$, where $s' = s + \frac{1}{2} - \frac{1}{\pi}$. Furthermore, for any $f, h \in V$, let
$$d(f, h) = \|f - h\|_{L^p(g)}$$
be the weighted $L^p$-norm on $V$. Then,
$$d(f_{jk}, f_{jk'}) = \beta_{j,k} 2^{j(1/2 - 1/p)} \|\psi\|_p =: \delta.$$

Plugging-in $f_0 \equiv 0$ and $f = f_{jk}$ in (6.1), we obtain

$$(6.2) \quad -2\ln \Lambda_n(f_0, f_{jk}) = \left(\frac{f_{jk}(\mathbf{X})}{\sigma(\mathbf{X})}\right)' \Xi^{-1} \left(\frac{f_{jk}(\mathbf{X})}{\sigma(\mathbf{X})}\right) + 2\left(\frac{f_{jk}(\mathbf{X})}{\sigma(\mathbf{X})}\right)' \Xi^{-1} \varepsilon.$$

Write
$$\ln \Lambda_n(f_0, f) = \{\ln \Lambda_n(f_0, f_{jk}) + \lambda_n\} - \lambda_n =: u_{nk} - v_{nk}.$$

Note, also, that the first component (6.2) is nonnegative, since $\Xi$ (and so $\Xi^{-1}$) is positive definite.

By the Cauchy–Schwarz inequality and Lemma 6.2 below, we obtain

$$(6.3) \quad \mathrm{E}\left[\left|\varepsilon' \Xi^{-1}\left(\frac{f_{jk}(\mathbf{X})}{\sigma(\mathbf{X})}\right)\right|\right] \leq \left\{\mathrm{E}\left[\left(\frac{f_{jk}(\mathbf{X})}{\sigma(\mathbf{X})}\right)' \Xi^{-1} \left(\frac{f_{jk}(\mathbf{X})}{\sigma(\mathbf{X})}\right)\right]\right\}^{1/2}.$$

Therefore, by (6.2), (6.3), Chebyshev inequality and the aforementioned positivity of the component in (6.2), we obtain

$$P(u_{nk} > 0) = P(\ln \Lambda_n(f_0, f) > -\lambda_n)$$
$$(6.4) \qquad \geq 1 - \lambda_n^{-1} \mathrm{E}\left[\frac{1}{2}\left(\frac{f_{jk}(\mathbf{X})}{\sigma(\mathbf{X})}\right)' \Xi^{-1}\left(\frac{f(\mathbf{X})}{\sigma(\mathbf{X})}\right) + \left|\varepsilon' \Xi^{-1}\left(\frac{f_{jk}(\mathbf{X})}{\sigma(\mathbf{X})}\right)\right|\right]$$
$$\geq 1 - \frac{A + \sqrt{2A}}{2\lambda_n}.$$

Now, $\mathbf{1}'\Xi\mathbf{1} = \sum_{i,l} \xi_{il} = \mathrm{Var}(\sum_{i=1}^n \varepsilon_i) \sim c_\alpha n^{2-\alpha}$ via (1.3). Also,
$$(\mathbf{1}'\Xi^{-1}\mathbf{1})(\mathbf{1}'\Xi\mathbf{1}) = n^2,$$

so that

(6.5) $$\mathbf{1}'\Xi^{-1}\mathbf{1} \sim c_\alpha^{-1} n^\alpha.$$

Furthermore,

(6.6)
$$\begin{aligned}
\mathrm{E}&\left[\left(\frac{f_{jk}(\mathbf{X})}{\sigma(\mathbf{X})}\right)'\Xi^{-1}\left(\frac{f_{jk}(\mathbf{X})}{\sigma(\mathbf{X})}\right)\right]\\
&= \mathrm{E}[f_{jk}^2(X)/\sigma^2(X)]\sum_{i=1}^n \xi_{ii}^{-1} + \{\mathrm{E}[f_{jk}(X)/\sigma(X)]\}^2 \sum_{i\neq l}\xi_{il}^{-1}\\
&\leq \|1/\sigma\|_\infty^2 \|\psi\|_2^2 \beta_{j,k}^2 \operatorname{trace}(\Xi^{-1}) + 2^{-j}\beta_{j,k}^2(\mathbf{1}'\Xi^{-1}\mathbf{1})\|1/\sigma\|_\infty^2\\
&= O(n)\beta_{j,k}^2 + O(2^{-j}n^\alpha)\beta_{j,k}^2 = O(n)\beta_{j,k}^2.
\end{aligned}$$

Summarizing, we obtain that the nominator in (6.4) is bounded by $C_1 n \beta_{j,k}^2$. We now choose $j$, according to

(6.7) $$2^j = C_2 \left(\frac{n}{\log n}\right)^{1/(2s')}.$$

Then,

$$j \ln 2 \geq \frac{1}{2(s+1/2-1/\pi)}(\log n - \log\log n) + \log C > \frac{\log n}{4(s+1/2-1/\pi)} =: \lambda_n.$$

Therefore,

$$P(u_{nk} > 0) > 1 - 4C_1 C_2^{-2s'}(s+1/2-1/\pi) > \pi_0 > 0$$

by the appropriate choice of $C_2$ in (6.7). Consequently,

(6.8)
$$\begin{aligned}
\inf_{\hat{f}_n} \sup_{f:\, f\circ G^{-1}\in\mathcal{B}_{\pi,r}^s} &\mathrm{E}\|f-\hat{f}_n\|_{L^p(g)}^p\\
&\geq \inf_{\hat{f}_n}\sup_{f:\,f\circ G^{-1}\in\mathcal{B}_{\pi,r}^s} P(\|f-\hat{f}_n\|_{L^p(g)} > \delta/2)\\
&\geq C_p \pi_0 \delta^p = C_p 2^{-jp(s-1/\pi+1/p)} = C_p\left(\frac{\log n}{n}\right)^{-p/2\alpha_S}.
\end{aligned}$$

*Dense case.* Let $\eta$ be the vector with components $\eta_k = \pm 1$, $k=0,\ldots,2^j-1$. Let $\eta^i$ be the vector with components $\eta_k^i = (-1)^{1\{i=k\}}\eta_k$. Let $f_{j\eta}(x) = \gamma_j \times \sum_{k=0}^{2^j-1}\eta_k \psi_{j,k}(G(x))$. To have $f_{j\eta}\circ G^{-1} \in \mathcal{B}_{\pi,r}^s$, we must have $\gamma_j \leq A 2^{-j(s+1/2)}$.

Note that $f_{j\eta} - f_{j\eta^i} = \pm \gamma_j \psi_{ji}$. Now, plug-in $f = f_{j\eta}$ and $f_0 = f_{j\eta^i}$ in (6.1) to get

$$-2\ln \Lambda_n(f_0,f) = \gamma_j^2\left(\frac{\psi_{ji}(\mathbf{X})}{\sigma(\mathbf{X})}\right)'\Xi^{-1}\psi_{ji}(\mathbf{X}) \pm 2\gamma_j\left(\frac{\psi_{ji}(\mathbf{X})}{\sigma(\mathbf{X})}\right)'\Xi^{-1}\varepsilon.$$

As in (6.6), we have

$$\mathrm{E}_{f_{j\eta}}[|\ln \Lambda_n(f_0, f)|] \le C_3 n \gamma_j^2 \le \pi_0,$$

if we choose

$$2^j \sim C_4 n^{1/(2s+1)}$$

with the appropriate $C_4$. Now, as in [25],

$$\inf_{\hat{f}_n} \max_\eta \mathrm{E}_{f_{j\eta}} \|\hat{f}_n - f_{j\eta}\|_{l^1(g)} \ge C 2^{j/2} \gamma_j,$$

which, by Cauchy–Schwarz inequality, yields

(6.9) $$\inf_{\hat{f}_n} \sup_{f \in \mathcal{B}_{\pi,r}^s \cap \{0\}} \mathrm{E}\|f - \hat{f}_n\|_{L^p(g)}^p \ge C_p n^{-p/2\alpha_D}.$$

Therefore, via (6.9) and (6.8), we obtain the i.i.d. lower bounds in Theorem 3.1. It finishes the proof in case of $f \in \mathcal{B}_{\pi,r}^s \cap \Delta$. If $f \notin \Delta$, then its mean has to be estimated. The lower bound follows in the very same way as on page 645 of [29]. This finishes the proof of Theorem 3.1.

LEMMA 6.2. *We have*

$$\mathrm{E}[|\varepsilon' \Xi^{-1} f(\mathbf{X})|^2] \le \mathrm{E}[f(\mathbf{X})' \Xi^{-1} f(\mathbf{X})].$$

PROOF. Bearing in mind the symmetry of $\Xi$,

$$\mathrm{E}[|\varepsilon' \Xi^{-1} f(\mathbf{X})|^2]$$
$$= \mathrm{E}\left[\sum_{i,l} \sum_{i_1,l_1} \varepsilon_i \varepsilon_{i_1} \xi_{il}^{-1} \xi_{i_1 l_1}^{-1} f(X_l) f(X_{l_1})\right]$$
$$= \mathrm{E}[f^2(X)] \sum_{i,l,i_1} \mathrm{E}[\varepsilon_i \varepsilon_{i_1}] \xi_{il}^{-1} \xi_{i_1 l}^{-1} + \{\mathrm{E}[f(X)]\}^2 \sum_{i,l,i_1 l_1 \ne l} \mathrm{E}[\varepsilon_i \varepsilon_{i_1}] \xi_{il}^{-1} \xi_{i_1 l_1}^{-1}$$
$$= \mathrm{E}[f^2(X)] \sum_{i_1,l} \xi_{i_1 l}^{-1} \sum_i \xi_{i_1 i} \xi_{il}^{-1} + \{\mathrm{E}[f(X)]\}^2 \sum_{l,i_1,l_1 \ne l} \xi_{i_1 l_1}^{-1} \sum_i \xi_{i_1 i} \xi_{il}^{-1}$$
$$= \mathrm{E}[f^2(X)] \sum_{i_1,l} \xi_{i_1 l}^{-1} (\Xi \Xi^{-1})_{i_1 l} + \{\mathrm{E}[f(X)]\}^2 \sum_{l,i_1,l_1 \ne l} \xi_{i_1 l_1} (\Xi \Xi^{-1})_{i_1 l}$$
$$= \mathrm{E}[f^2(X)] \operatorname{trace}(\Xi^{-1}) + \{\mathrm{E}[f(X)]\}^2 \sum_{l,l_1 \ne l} \xi_{l l_1}^{-1}$$
$$= \mathrm{E}[f(\mathbf{X})' \Xi^{-1} f(\mathbf{X})]. \qquad \square$$

## 7. Proofs: upper bounds.

7.1. *Decomposition of empirical wavelet coefficients.* Here, we establish decomposition of the form,

$$\hat{\beta}_{j,k} - \beta_{j,k} = \text{i.i.d. part} + \text{martingale part} + \text{wavelet LRD part}.$$

From (4.3),

(7.1) $\quad E[\hat{\beta}_{j,k}] = E[\psi_{j,k}(X_1)f(X_1)] = \int \psi_{j,k}(y) f(G^{-1}(y))\,dy = \beta_{j,k}.$

We set $U_i := G(X_i), i = 1,\ldots,n$, the $U_i$'s are uniformly distributed on $[0,1]$, by independence

$$E[\psi_{j,k}(U_1)\sigma(X_1)\varepsilon_1] = E[\psi_{j,k}(U_1)\sigma(X_1)]E[\varepsilon_1] = 0,$$

(7.2)
$$\begin{aligned}
\hat{\beta}_{j,k} - \beta_{j,k} &= \frac{1}{n}\sum_{i=1}^{n}(\psi_{j,k}(U_i)Y_i - E[\psi_{j,k}(U_i)Y_i]) \\
&= \frac{1}{n}\sum_{i=1}^{n}(\psi_{j,k}(U_i)f(X_i) - E[\psi_{j,k}(U_1)f(X_1)]) \\
&\quad + \frac{1}{n}\sum_{i=1}^{n}\psi_{j,k}(U_i)\sigma(X_i)\varepsilon_i \\
&=: A_0 + A_1.
\end{aligned}$$

Note that $A_0$ is the sum of i.i.d. random variables, whereas the dependence structure is included in $A_1$ only. The part $A_1$ is decomposed further. Let $\mathcal{F}_i = \sigma(\eta_i, X_i, \eta_{i-1}, X_{i-1}, \ldots)$. Let $\varepsilon_{i,i-1} = \varepsilon_i - \eta_i$. Note that $\varepsilon_{i,i-1}$ is $\mathcal{F}_{i-1}$-measurable and $(\eta_i, X_i)$ is independent of $\mathcal{F}_{i-1}$. Thus,

$$E[\psi_{j,k}(U_i)\sigma(X_i)\varepsilon_i | \mathcal{F}_{i-1}] = \varepsilon_{i,i-1} E[\psi_{j,k}(U_1)\sigma(X_1)].$$

We write

(7.3)
$$\begin{aligned}
A_1 &= \frac{1}{n}\sum_{i=1}^{n}\psi_{j,k}(U_i)\sigma(X_i)\varepsilon_i \\
&= \frac{1}{n}\sum_{i=1}^{n}(\psi_{j,k}(U_i)\sigma(X_i)\varepsilon_i - E[\psi_{j,k}(U_i)\sigma(X_i)\varepsilon_i | \mathcal{F}_{i-1}]) \\
&\quad + \frac{1}{n}E[\psi_{j,k}(U_1)\sigma(X_1)]\sum_{i=1}^{n}\varepsilon_{i,i-1} \\
&=: A_2 + A_3
\end{aligned}$$

and

$$\hat{\beta}_{j,k} - \beta_{j,k} = A_0 + A_2 + A_3$$
(7.4)
$$=: \text{i.i.d. part} + \text{martingale part} + \text{wavelet LRD part}.$$

Consider, also, the following corresponding decomposition for the scaling coefficients $\alpha_{j,k}$:

$$\hat{\alpha}_{j,k} - \alpha_{j,k} = B_0 + B_2 + B_3$$
(7.5)
$$=: \text{i.i.d. part} + \text{martingale part} + \text{scaling LRD part}.$$

An important feature of this decomposition is that the LRD term involves the partial sums of $\varepsilon_{i,i-1}$ only. Furthermore, if (1.10) holds, then $A_3 \equiv 0$ and the LRD part does not contribute. On the other hand, the scaling LRD part is always present.

As for the shape estimation, let $\alpha^*_{j,k}$ be the scaling coefficient of $f^* \circ G^{-1}$. Clearly,

$$\alpha^*_{j,k} = \alpha_{j,k} - \int_0^1 f(G^{-1}(y))\,dy \int_0^1 \phi_{j,k}(y)\,dy =: \alpha_{j,k} - c_{f,G}\mathrm{E}[\phi_{j,k}(U_1)].$$

Let $\hat{c}_{f,G}$ be an estimator of $c_{f,G}$ [e.g., $\hat{c}_{f,G} = \frac{1}{n}\sum_{i=1}^n f(X_i)$]. Then, we decompose

$$\hat{\alpha}^*_{j,k} - \alpha^*_{j,k} = B_0 + B_2 + B_3$$
$$= \text{i.i.d. part} + \text{martingale part} + c_{f,G}\mathrm{E}\phi_{j,k}(U_1) - \hat{c}_{f,G}\mathrm{E}\phi_{j,k}(U_1)$$
$$+ \frac{1}{n}\mathrm{E}[\phi_{j,k}(U_1)\sigma(X_1)]\sum_{i=1}^n \varepsilon_{i,i-1} - \frac{1}{n}\mathrm{E}[\phi_{j,k}(U_1)]\sum_{i=1}^n \sigma(X_i)\varepsilon_i.$$

If $\sigma(\cdot) \equiv 1$, then the last two terms equal

$$\frac{1}{n}\mathrm{E}[\phi_{j,k}(U_1)]\sum_{i=1}^n \varepsilon_{i,i-1} - \frac{1}{n}\mathrm{E}[\phi_{j,k}(U_1)]\sum_{i=1}^n \varepsilon_i = -\frac{1}{n}\mathrm{E}[\phi_{j,k}(U_1)]\sum_{i=1}^n \eta_i,$$

which is the just sum of i.i.d. random variables. Consequently, if (1.10) holds, then the LRD effect is not present in the scaling coefficient estimation. Otherwise, the LRD part is present and affects convergence rates. Therefore, by removing the scaling coefficient $\phi_{0,0}$, we guarantee that LRD does not affect the shape estimation.

7.2. *Decomposition of the modified wavelet coefficients.* In this section, we decompose $\tilde{\beta}_{j,k}$. Let us redefine

$$\mathcal{F}_i = \sigma(\eta_i, X_i, \eta_{i-1}, X_{i-1}, \ldots) \vee \sigma(X'_1, \ldots, X'_n).$$

Note that
$$\mathrm{E}[\psi_{j,k}(\hat{G}_n(X_i))\sigma(X_i)\varepsilon_i|\mathcal{F}_{i-1}] = \mathrm{E}[\psi_{j,k}(\hat{G}_n(X_1))\sigma(X_1)]\varepsilon_{i,i-1}$$
and $\psi_{j,k}(\hat{G}_n(X_i))\sigma(X_i)\varepsilon_i$ is $\mathcal{F}_i$-measurable. [This shows the importance of defining $\hat{G}_n(\cdot)$ based on the first different of the sample, $X'_1, \ldots, X'_n$.] Similarly to (7.2) and (7.3), we decompose

$$
\begin{aligned}
\tilde{\beta}_{j,k} - \beta_{j,k} &= \frac{1}{n}\sum_{i=1}^{n}(\psi_{j,k}(\hat{G}_n(X_i))Y_i - \beta_{j,k}) \\
&= \frac{1}{n}\sum_{i=1}^{n}(\psi_{j,k}(\hat{G}_n(X_i))f(X_i) - \beta_{j,k}) \\
&\quad + \frac{1}{n}\sum_{i=1}^{n}(\psi_{j,k}(\hat{G}_n(X_i))\sigma(X_i)\varepsilon_i \\
&\quad\quad - \mathrm{E}[\psi_{j,k}(\hat{G}_n(X_i))\sigma(X_i)\varepsilon_i|\mathcal{F}_{i-1}]) \\
&\quad + \frac{1}{n}\mathrm{E}[\psi_{j,k}(\hat{G}_n(X_1))\sigma(X_1)]\sum_{i=1}^{n}\varepsilon_{i,i-1} \\
&=: \tilde{A}_0 + \tilde{A}_2 + \tilde{A}_3.
\end{aligned}
$$
(7.6)

7.3. *Moment bounds.*

LEMMA 7.1. *For all $j \geq 0$ and $k = 0, \ldots, 2^j - 1$ and $p \geq 2$,*
$$\mathrm{E}[|\hat{\beta}_{j,k} - \beta_{j,k}|^p] = O(n^{-p/2}) + O(2^{-jp/2}n^{-p\alpha/2}) \tag{7.7}$$
*as long as $2^j \leq n$. The bound also applies to scaling coefficients $|\hat{\alpha}_{j,k} - \alpha_{j,k}|^p$. Moreover, if (1.10) holds, then*
$$\mathrm{E}[|\hat{\beta}_{j,k} - \beta_{j,k}|^p] = O(n^{-p/2}).$$

PROOF. *I.i.d. part.* By using Rosenthal's inequality, [16], page 132,
$$
\begin{aligned}
\mathrm{E}|A_0|^p &= \frac{1}{n}\mathrm{E}\left|\sum_{i=1}^{n}(\psi_{j,k}(G(X_i))f(X_i) - \mathrm{E}[\psi_{j,k}(G(X_i))f(X_i)])\right|^p \\
&\leq Cn^{-p}\|f\|_\infty^p(n2^{j(p/2-1)} + n^{p/2}) = O(n^{-p/2})
\end{aligned}
$$
(7.8)

as long as $2^j \leq n$.

*LRD part.* If (1.10), then the LRD part vanishes. Otherwise, note that
$$\mathrm{E}[|\psi_{j,k}(U_1)\sigma(X_1)|^p] \leq \|\sigma\|_\infty^p\|\psi\|_p^p 2^{j(p/2-1)}. \tag{7.9}$$

Since $\sum_{i=1}^{n} \varepsilon_{i,i-1}$ is a centered normal random variable with variance

$$(7.10) \quad v_n^2 := \operatorname{Var}\left(\sum_{i=1}^{n} \varepsilon_{i,i-1}\right) \sim d_\alpha n^{2-\alpha},$$

we obtain

$$(7.11) \quad \mathrm{E}|A_3|^p = O(2^{-jp/2} n^{-p\alpha/2}).$$

*Martingale part.* In the light of the decomposition (7.4), we see that $nA_2 =: \sum_{i=1}^{n} d_i$ is a martingale, where

$$d_i = \psi_{j,k}(U_i)\sigma(X_i)\varepsilon_i - \mathrm{E}[\psi_{j,k}(U_i)\sigma(X_i)\varepsilon_i|\mathcal{F}_{i-1}]$$
$$= \varepsilon_{i,i-1}(\psi_{j,k}(U_i)\sigma(X_i) - \mathrm{E}[\psi_{j,k}(U_1)\sigma(X_1)]) + \eta_i \psi_{j,k}(U_i)\sigma(X_i).$$

Note that the first and the second term are uncorrelated, both unconditionally and conditionally on $\mathcal{F}_{i-1}$. By (7.9),

$$\mathrm{E}[|d_i|^p] \leq 2^{p-1}(\mathrm{E}[|\varepsilon_{i,i-1}|^p]\mathrm{E}[|\psi_{j,k}(U_i)\sigma(X_i) - \mathrm{E}[\psi_{j,k}(U_i)\sigma(X_i)]|^p]$$
$$+ \mathrm{E}[|\eta_i|^p]\mathrm{E}[|\psi_{j,k}(U_i)\sigma(X_i)|^p])$$
$$\leq C\mathrm{E}[|\psi_{j,k}(U_1)\sigma(X_1)|^p] = C2^{j(p/2-1)}.$$

Now,

$$(7.12) \quad \begin{aligned} \sigma_i^2 &:= \mathrm{E}[d_i^2|\mathcal{F}_{i-1}] \\ &= \mathrm{E}[\psi_{j,k}^2(U_1)\sigma^2(X_1)]\mathrm{E}[\eta_1^2] + \varepsilon_{i,i-1}^2 \operatorname{Var}[\psi_{j,k}(U_1)\sigma(X_1)]. \end{aligned}$$

Using $\mathrm{E}[\psi_{j,k}^2(U_1)\sigma(X_1)] = O(1)$,

$$\mathrm{E}\left[\left(\sum_{i=1}^{n} \mathrm{E}(d_i^2|\mathcal{F}_{i-1})\right)^{p/2}\right]$$
$$= \mathrm{E}\left(n\mathrm{E}[\psi_{j,k}^2(U_1)\sigma(X_1)]\mathrm{E}\eta_1^2 + \operatorname{Var}[\psi_{j,k}(U_1)\sigma(X_1)]\sum_{i=1}^{n} \varepsilon_{i,i-1}^2\right)^{p/2}$$
$$\leq C_p n^{p/2}(\mathrm{E}\psi_{j,k}^2(U_1))^{p/2} + C_p(\operatorname{Var}\psi_{j,k}(U_1))^{p/2}\mathrm{E}\left(\sum_{i=1}^{n} \varepsilon_{i,i-1}^2\right)^{p/2}$$
$$\leq C n^{p/2}.$$

Using Rosenthal's inequality for martingales [14], page 25,

$$(7.13) \quad \begin{aligned} \mathrm{E}|A_2|^p &\leq Cn^{-p}\mathrm{E}\left(\sum_{i=1}^{n} \mathrm{E}(d_i^2|\mathcal{F}_{i-1})\right)^{p/2} + Cn^{-p}\sum_{i=1}^{n} \mathrm{E}|d_i|^p \\ &\leq C(n^{-p/2} + n^{-p} n 2^{j(p/2-1)}) \leq Cn^{-p/2} \end{aligned}$$

as long as $2^j \leq n$. Now, (7.7) follows from (7.8), (7.11) and (7.13). □

7.4. *Large deviation estimates.*

PROPOSITION 7.2. *Let $\lambda_{n,j}$ be as in (4.4). Assume that $j$ is such that $2^j < (n/\log n)$. For any $r > 0$, there exist positive constants $\tau$ and $C(r, p, \tau)$ such that*

(7.14) $$P(|\hat{\beta}_{j,k} - \beta_{j,k}| > \tau \lambda_{n,j}/2) \leq C(r, p, \tau) n^{-rp}.$$

*A similar bound is valid for $\hat{\alpha}_{j,k} - \alpha_{j,k}$.*

PROOF. We obtain (7.14) separately for $A_0, A_3$ and $A_2$ and apply triangular inequalities in (7.4) to complete the proof. A similar approach works for (7.5).

*I.i.d. part.* For $A_0$, we have from the Bernstein inequality as long as $2^j \leq (n/\log n)$ (see, e.g., [19], Proposition 3)

(7.15) $$P\left(|A_0| > \frac{\tau}{2}\sqrt{\frac{\log n}{n}}\right) \leq 2\exp\left(-\frac{3\tau^2 \log n}{8\|f\|_\infty \max\{3, \tau\}}\right)$$

for all $n$. The bound in (7.14) is valid for the i.i.d. part with

(7.16) $$\tau \geq \max\{\tfrac{8}{3}\|f\|_\infty rp, \sqrt{8rp\|f\|_\infty}\}.$$

*LRD part.* First, if (1.10) holds, then LRD part vanishes. Otherwise, we recall (7.9) and that $\sum_{i=1}^n \varepsilon_{i,i-1}$ is a centered normal r.v. with variance (7.10). For sufficiently large $n$,

$$P\left(|A_3| > \tau\sqrt{\frac{\log n}{n}}/2\right) \leq C\exp\left(-\frac{\tau n \log n 2^j}{4d_\alpha \|\sigma\|_\infty^2 \|\psi\|_1^2 n^{2-\alpha}}\right).$$

Therefore, for all $j$ such that $2^j > n^{1-\alpha}$,

$$P\left(|A_3| > \tau\sqrt{\frac{\log n}{n}}/2\right) \leq C\exp\left(-\frac{\tau \log n}{4d_\alpha \|\sigma\|_\infty^2 \|\psi\|_1^2}\right) \leq Cn^{-rp}$$

for all $n$, if

(7.17) $$\tau > 4d_\alpha rp\|\sigma\|_\infty^2 \|\psi\|_1^2.$$

If, now, $2^j < n^{1-\alpha}$, then

$$P(|A_3| > \tau \tilde{\lambda}_{n,j}/2) = P\left(\left|\sum_{i=1}^n \varepsilon_{i,i-1}\right| > \frac{\tau n \tilde{\lambda}_{n,j}}{2|\mathrm{E}[\psi_{j,k}(U_1)\sigma(X_1)]|}\right)$$
$$\leq C\exp\left(-\frac{\tau \log n}{4d_\alpha \|\sigma\|_\infty^2 \|\psi\|_1^2}\right) \leq Cn^{-rp}$$

for the same choice of $\tau$ as in (7.17).

*Martingale part.* For $A_2$, we will use a new Bernstein's inequality for martingales. We recall the following lemma from [11].

LEMMA 7.3. *Let $(d_i, \mathcal{F}_i), i \geq 1$, be a martingale difference sequence. Denote $\sigma_i^2 = \mathrm{E}[d_i^2 | \mathcal{F}_{i-1}]$. For any $x, L, a > 0$,*

$$P\left(\left|\sum_{i=1}^n d_i\right| > x, \sum_{i=1}^n d_i^2 \mathbb{I}_{\{|d_i| > a\}} + \sum_{i=1}^n \sigma_i^2 \leq L\right) \leq 2\exp\left(-\frac{x^2}{2(L + ax/3)}\right).$$

We apply this lemma to our martingale sequence $d_i$ and $\sigma_i^2$ defined in (7.12), with a very precise choice of truncation levels $a$ and $L$ (clearly, they cannot be too big). Let

$$H_n = H_n(a) := \sum_{i=1}^n d_i^2 \mathbb{I}_{\{|d_i| > a\}} + \sum_{i=1}^n \sigma_i^2,$$

(7.18)
$$P\left(\frac{1}{n}\left|\sum_{i=1}^n d_i\right| > x\right) \leq P\left(\frac{1}{n}\left|\sum_{i=1}^n d_i\right| > x, H_n \leq L\right) + P(H_n \geq L)$$

$$\leq 2\exp\left(-\frac{n^2 x^2}{2(L + anx/3)}\right) + P(H_n > L).$$

We take

(7.19)  $L := L_n = 2n(A\log n + \mathrm{E}[\psi_{j,k}^2(U_1)\sigma^2(X_1)]\mathrm{E}[\eta_1^2]) =: 2n(A\log n + C_1)$

with $A > 0$ to be specified below,

$$P\left(\sum_{i=1}^n \sigma_i^2 > L/2\right) = P\bigg(\mathrm{Var}[\psi_{j,k}(U_1)\sigma(X_1)] \sum_{i=1}^n \varepsilon_{i,i-1}^2$$

$$+ \mathrm{E}[\psi_{j,k}^2(U_1)\sigma^2(X_1)]\mathrm{E}[\eta_1^2]n > L/2\bigg)$$

$$= P\bigg(\mathrm{Var}[\psi_{j,k}(U_1)\sigma(X_1)] \sum_{i=1}^n \varepsilon_{i,i-1}^2 > An\log n\bigg)$$

(7.20)
$$\leq P\bigg(\bigcup_{i=1}^n \Big\{\varepsilon_{i,i-1}^2 > \frac{A\log n}{\mathrm{Var}[\psi_{j,k}(U_1)\sigma(X_1)]}\Big\}\bigg)$$

$$\leq nP\bigg(\varepsilon_{1,0}^2 > \frac{A\log n}{\mathrm{Var}[\psi_{j,k}(U_1)\sigma(X_1)]}\bigg)$$

$$\leq Cn\exp\bigg(-\frac{A\log n}{2\mathrm{Var}[\varepsilon_{1,0}]\mathrm{Var}[\psi_{j,k}(U_1)\sigma(X_1)]}\bigg)$$

$$= Cn^{-rp}$$

by the choice

(7.21) $$A = 2(rp+1)\operatorname{Var}[\varepsilon_{1,0}]\operatorname{Var}[\psi_{j,k}(U_1)\sigma(X_1)].$$

Further, note that

$$d_i^2 \leq 4((\psi_{j,k}(U_i)\sigma(X_i) - \mathrm{E}[\psi_{j,k}(U_1)\sigma(X_1)])^2 \varepsilon_{i,i-1}^2 + \eta_i^2 \psi_{j,k}^2(U_i)\sigma^2(X_i)).$$

Thus, for any $a > (A\log n)^{1/2}$,

(7.22)
$$\begin{aligned}
P\left(\sum_{i=1}^n d_i^2 \mathbb{I}_{\{|d_i|>a\}} > L/2\right) &\leq P\left(\sum_{i=1}^n d_i^2 \mathbb{I}_{\{|d_i|>a\}} > An\log n\right)\\
&\leq nP(d_1^2 \mathbb{I}_{\{|d_1|>a\}} > A\log n)\\
&\leq nP(d_1^2 > (A\log n) \vee a^2)\\
&\leq nP(d_1^2 > a^2) = nP(d_1^2 > a^2)\\
&\leq nP((\psi_{j,k}(U_1)\sigma(X_1)\\
&\quad - \mathrm{E}[\psi_{j,k}(U_1)\sigma(X_1)])^2 \varepsilon_{1,0}^2 > a^2/2)\\
&\quad + nP(\eta_1^2 \psi_{j,k}^2(U_1)\sigma^2(X_1) > a^2/2).
\end{aligned}$$

Since

$$|\psi_{j,k}^2(U_1)\sigma^2(X_1)| \leq 2^j \|\psi\|_\infty^2 \|\sigma\|_\infty^2 =: C_0 2^j \leq C_0 \frac{n}{\log n},$$

we have

(7.23)
$$\begin{aligned}
P(\eta_1^2 \psi_{j,k}^2(U_1)\sigma^2(X_1) > a^2/2) &\leq P\left(\eta_i^2 > \frac{a^2}{2C_0 2^j}\right)\\
&\leq C\exp\left(-\frac{a^2 \log n}{4C_0 n}\right) \leq Cn^{-rp},
\end{aligned}$$

by choosing $a = B\sqrt{n}$ with

(7.24) $$B = 4C_0 rp.$$

A similar bound applies to the first term in (7.22).

Combining (7.18), (7.20) and (7.23),

(7.25) $$P\left(\frac{1}{n}\left|\sum_{i=1}^n d_i\right| > x\right) \leq 2\exp\left(-\frac{n^2 x^2}{2(L+anx/3)}\right) + Cn^{-rp},$$

where $L$ as in (7.19). Take, now, $x = \frac{\tau}{2}\frac{\log n}{\sqrt{n}}$, and note that

$$\frac{n^2 x^2}{2(L+anx/3)} \leq \frac{n(\log n)^2 \tau^2/8}{n\log n(A+C_1) + B\tau/6 n\log n},$$

so that (7.14) follows for the martingale part by taking

$$\tau \geq \max\{\sqrt{8(A+C_1)rp}, Brp\tfrac{8}{6}\},$$

where $A$, $C_1$ and $B$ were defined in (7.21), (7.19) and (7.24), respectively. □

7.5. *Bounds for the modified wavelet coefficients.* Let us start with the following bound:

LEMMA 7.4. *For all $2^j \leq \sqrt{n}$, we have*

(7.26) $\quad \mathrm{E}[|\psi_{j,k}(\hat{G}_n(X_1))\sigma(X_1)|^p] \leq C_p \|\sigma\|_\infty^p \|\psi\|_p^p 2^{j(p/2-1)},$

*where $C_p$ is a constant depending only on $p$.*

PROOF. Let $\Theta_n(x)$ be a random element between $\hat{G}_n(x)$ and $G(x)$. Then,

$$\begin{aligned}
\mathrm{E}[|(\psi_{j,k}(\hat{G}_n(X_1)) &- \psi_{j,k}(G(X_i)))\sigma(X_1)|^p] \\
&\leq \mathrm{E}\mathrm{E}[|\psi'_{j,k}(G(X_1))|^p |\Theta_n(X_1)|^p |\sigma(X_1)|^p | X_1] \\
&\leq \|\sigma\|_\infty^p \mathrm{E}\Big[\sup_x |\Theta_n(x)|^p\Big] \mathrm{E}[|\psi'_{j,k}(G(X_1))|^p] = O(n^{-p/2} 2^{3jp/2 - j}).
\end{aligned}$$

In the above computation, we used independence of $\hat{G}_n(\cdot)$ of $X_1$ and the standard bound on the supremum norm of the empirical process. Consequently,

$$\mathrm{E}[|\psi_{j,k}(\hat{G}_n(X_1))\sigma(X_1)|^p] \leq \mathrm{E}[|\psi_{j,k}(G(X_1))\sigma(X_1)|^p] + O(n^{-p/2} 2^{3jp/2 - j})$$

and the bound is of order $2^{j(p/2-1)}$ if and only if $2^j \leq \sqrt{n}$. □

With help of the above lemma, we conclude that the results for $\hat{\beta}_{j,k}$ can be rewritten for $\tilde{\beta}_{j,k}$.

LEMMA 7.5. *Assume that $\|f \circ G^{-1}\|_{\mathrm{Lip}(1/2)} < \infty$. The bounds of Lemma 7.1 and Proposition 7.2 remain valid for $\tilde{\beta}_{j,k}$ and $\tilde{\alpha}_{j,k}$ as long as $2^j \leq \sqrt{n}$.*

PROOF. The bounds for the first part of the decomposition (7.6), $\tilde{A}_0$, follow from [19], Proposition 6. To deal with the LRD part, $\tilde{A}_3$, we simply replace (7.9) with (7.26) [see (7.11) and the computation leading to (7.17)]. Similarly, note that the moment bounds and large deviations for the martingale part involve only the behavior of $\mathrm{E}[|\psi_{j,k}(\hat{G}_n(X_1))\sigma(X_1)|^p]$ instead of $\mathrm{E}[|\psi_{j,k}(G(X_1))\sigma(X_1)|^p]$. □

7.6. *Proof of Theorem 4.1.* In what follows, $D_j = \{k, k = 0, 1, \ldots, 2^j - 1\}$, we split $\hat{f}_n - f$ into three parts,

$$\mathrm{E}\|f - \hat{f}_n\|_{L^p(g)}^p$$
$$\leq 3^{p-1}\Bigg(\mathrm{E}\bigg\|\sum_{k \in D_{j_0}} (\alpha_{j_0,k} - \hat{\alpha}_{j_0,k})\phi_{j_0,k}(G(\cdot))\bigg\|_{L^p(g)}^p$$
$$+ \mathrm{E}\bigg\|\sum_{j=j_0}^{j_1}\sum_{k \in D_j}\beta_{j,k}\psi_{j,k}(G(\cdot))$$
$$- \sum_{j=j_0}^{j_1}\sum_{k \in D_j}\hat{\beta}_{j,k}\mathbb{I}_{\{|\hat{\beta}_{j,k}|>\tau_0\lambda_{n,j}\}}\psi_{j,k}(G(\cdot))\bigg\|_{L^p(g)}^p$$
$$+ \bigg\|\sum_{j \geq j_1}\sum_{k \in D_j}\beta_{j,k}\psi_{j,k}(G(\cdot))\bigg\|_{L^p(g)}^p\Bigg)$$
$$:= \text{linear term} + \text{nonlinear term} + \text{bias term}.$$

*Bias term.* We use standard approximation results (see, e.g., [16], pages 123–124), introducing

(7.27) $$\delta := s - \left(\frac{1}{\pi} - \frac{1}{p}\right)_+ = s - \max\left(\frac{1}{\pi} - \frac{1}{p}, 0\right),$$

if $p \leq \pi$, $\delta = s$ and $\mathcal{B}_{\pi,r}^s \subseteq \mathcal{B}_{p,r}^s$, if $\pi < p$, $\delta = s - (\frac{1}{\pi} - \frac{1}{p})$ and $\mathcal{B}_{\pi,r}^s \subseteq \mathcal{B}_{p,r}^\delta$,

$$\bigg\|\sum_{j > j_1}\sum_{k \in D_j}\beta_{j,k}\psi_{j,k}(G(\cdot))\bigg\|_{L^p(g)}^p$$
$$= \int_0^1 \bigg|\sum_{j > j_1}\sum_{k \in D_j}\beta_{j,k}\psi_{j,k}(G(x))\bigg|^p g(x)\,dx$$
(7.28)
$$= \int_0^1 \bigg|\sum_{j > j_1}\sum_{k \in D_j}\beta_{j,k}\psi_{j,k}(u)\bigg|^p du$$
$$\leq C\|f \circ G^{-1}\|_{\mathcal{B}_{p,r}^\delta}^p 2^{-j_1\delta p} = O((\log n/n)^{\delta p}),$$

where we have used the definition (4.1) of $j_1$ for the last bound.

*The linear part.* Applying Lemma 7.1, the term $\mathrm{E}|\hat{\alpha}_{j_0,k} - \alpha_{j_0,k}|^p$ is proportional to $n^{-p\alpha/2}$. Therefore,

$$\mathrm{E}\bigg\|\sum_k (\alpha_{j_0,k} - \hat{\alpha}_{j_0,k})\phi_{j_0,k}(G(\cdot))\bigg\|_{L^p(g)}^p$$

$$\leq 2^{j_0(p/2-1)}\|\phi\|_p^p \sum_{k\in D_{j_0}} \mathrm{E}|\alpha_{j_0,k}-\alpha_{j_0,k}|^p$$

$$\leq C 2^{j_0(p/2-1)} 2^{j_0} \mathrm{E}|\hat\alpha_{j_0,k}-\alpha_{j_0,k}|^p = O(n^{-p\alpha/2}).$$

*Nonlinear term.* We follow the proof of Theorem 5.1 in [18], incorporating our moments and large deviations bounds accordingly. We refer to Appendix for the definition $l_{q,\infty}$ spaces. We use Temlyakov's property and Minkowski's inequality repeatedly.

$$\mathrm{E}\Big\|\sum_{(j,k)\in\Lambda_1}\beta_{j,k}\psi_{j,k}(G(\cdot))-\sum_{j,k\in\Lambda_1}\hat\beta_{(j,k)}\mathbb{I}_{\{|\hat\beta_{j,k}|>\tau_0\lambda_{n,j}\}}\psi_{j,k}(G(\cdot))\Big\|_{L^p(g)}^p$$

$$\leq 2^{p-1}\bigg(\mathrm{E}\Big\|\sum_{(j,k)\in\Lambda_1}(\beta_{j,k}-\hat\beta_{j,k})\mathbb{I}_{\{|\hat\beta_{j,k}|>\tau_0\lambda_{n,j}\}}\psi_{j,k}(G(\cdot))\Big\|_{L^p(g)}^p$$

$$+\mathrm{E}\Big\|\sum_{(j,k)\in\Lambda_1}\beta_{j,k}\mathbb{I}_{\{|\hat\beta_{j,k}|\leq\tau_0\lambda_{n,j}\}}\psi_{j,k}(G(\cdot))\Big\|_{L^p(g)}^p\bigg)$$

$$=: A+B.$$

Let us introduce some notation. We define $j_2$ to be such that $2^{j_2} = n^{1-\alpha}$. Further, let

$$\Lambda_2 = \{(j,k), j_2 \leq j \leq j_1, k=0,1,\ldots,2^j-1\}, \qquad \Lambda_3 = \Lambda_1 \setminus \Lambda_2.$$

*We start by the A-term.* Changing variables $u=G(x)$ we get

$$A \leq \mathrm{E}\int\bigg(\sum_{(j,k)\in\Lambda_1}|\hat\beta_{j,k}-\beta_{j,k}|^2\mathbb{I}_{\{|\hat\beta_{j,k}-\beta_{j,k}|>\tau\lambda_{n,j}/2\}}\psi_{j,k}^2(u)\bigg)^{p/2}du$$

$$+\mathrm{E}\int\bigg(\sum_{(j,k)\in\Lambda_1}|\hat\beta_{j,k}-\beta_{j,k}|^2\mathbb{I}_{\{|\beta_{j,k}|>\tau\lambda_{n,j}/2\}}\psi_{j,k}^2(u)\bigg)^{p/2}du$$

$$\leq \int\bigg\{\sum_{(j,k)\in\Lambda_1}[(\mathrm{E}|\hat\beta_{j,k}-\beta_{j,k}|^{2p}$$

$$\times P(|\hat\beta_{j,k}-\beta_{j,k}|>\tau\lambda_{n,j}/2))^{1/2}|\psi_{j,k}(u)|^p]^{2/p}\bigg\}^{p/2}du$$

$$+\int\bigg\{\sum_{(j,k)\in\Lambda_1}\mathbb{I}_{\{|\beta_{j,k}|>\tau\lambda_{n,j}/2\}}[\mathrm{E}|\hat\beta_{j,k}-\beta_{j,k}|^p]^{2/p}\psi_{j,k}^2(u)\bigg\}^{p/2}du$$

$$=: A_1+A_2.$$

Using the bounds of Lemma 7.1 and (A.2) below,

$$A_2 \le C \int \left\{ \sum_{(j,k)\in\Lambda_3} \mathbb{I}_{\{|\beta_{j,k}|>\tau_0\lambda_{n,j}/2\}} (2^{-jp/2} n^{-p\alpha/2})^{2/p} \psi_{j,k}^2(u) \right\}^{p/2} du$$

$$+ C \int \left\{ \sum_{(j,k)\in\Lambda_2} \mathbb{I}_{\{|\beta_{j,k}|>\tau_0\lambda_{n,j}/2\}} (n^{-p/2})^{2/p} \psi_{j,k}^2(u) \right\}^{p/2} du$$

$$\le C n^{-p\alpha/2} \sum_{j=0}^{j_2} 2^{j(p/2-1)} 2^{-jp/2} |D_j| \|\psi\|_p^p$$

$$+ Cc_n^p \sum_{j=j_2}^{j_1} \|\psi_{j,k}\|_p^p \sum_{k\in D_j} \mathbb{I}_{\{|\beta_{j,k}|>\tau_0\lambda_{n,j}/2\}}$$

$$\le C n^{-p\alpha/2} j_2 + C c_n^{p-q} \sup_{\lambda>0} \lambda^q \sum_{j=1}^{j_1} \sum_{k\in D_j} \|\psi_{j,k}\|_p^p \mathbb{I}_{\{|\beta_{j,k}|>\tau\lambda/2\}}$$

$$\le C n^{-p\alpha/2} \log n + C \tilde{\lambda}_n^{p-q} \|f\|_{l_{q,\infty}}^q.$$

In the second to last inequality, we used $c_n^q \le \tilde{\lambda}_n^q$ and the fact that for $j > j_2$ we have $\lambda_{n,j} = \tilde{\lambda}_n$.

As for $A_1$, we split this into 2 parts, according to $\Lambda_2$ and $\Lambda_3$. On $\Lambda_2$, using Lemma 7.1 and Proposition 7.2, we get (recall that then $\lambda_{n,j} = \tilde{\lambda}_n$)

$$(\mathrm{E}|\hat{\beta}_{j,k} - \beta_{j,k}|^{2p} P(|\hat{\beta}_{j,k} - \beta_{j,k}| > \tau_0 \lambda_{n,j}/2))^{1/2} = O((c_n^{2p} \tilde{\lambda}_n^{2p})^{1/2}) = \tilde{\lambda}_n^{2p}.$$

On $\Lambda_3$, we have

$$(\mathrm{E}|\hat{\beta}_{j,k} - \beta_{j,k}|^{2p} P(|\hat{\beta}_{j,k} - \beta_{j,k}| > \tau_0 \lambda_{n,j}/2))^{1/2}$$
$$= O(2^{-jp} n^{-\alpha p} (\log n)^{p/2}),$$

so that

$$A_1 \le C n^{-\alpha p} (\log n)^{p/2} \sum_{(j,k)\in\Lambda_3} 2^{-jp} \|\psi\|_{j,k}^p$$

(7.29) $$\qquad + C \tilde{\lambda}_n^{2p} \int \sum_{(j,k)\in\Lambda_2} |\psi_{j,k}(u)|^p\, du$$

$$\le C n^{-\alpha p} (\log n)^{p/2} 2^{-j_2 p/2} + C \tilde{\lambda}_n^{2p} \sum_{j=j_2}^{j_1} 2^j 2^{j(p/2-1)} = O(\tilde{\lambda}_n^p).$$

*For the B-term.*

$$B \leq \int \left\{ \sum_{(j,k)\in\Lambda_1} \mathbb{I}_{\{|\beta_{j,k}|>2\tau_0\lambda_{n,j}\}} P(|\hat{\beta}_{j,k}-\beta_{j,k}|>\tau_0\lambda_{n,j}/2)^{2/p} \beta_{j,k}^2 \psi_{j,k}^2(u) \right\}^{p/2} du$$

$$+ \int \left\{ \sum_{(j,k)\in\Lambda_1} \mathbb{I}_{\{|\beta_{j,k}|\leq 2\tau_0\lambda_{n,j}\}} \beta_{j,k}^2 \psi_{j,k}^2(u) \right\}^{p/2} du =: B_1 + B_2$$

and both terms are treated in the similar way as $A_1$ and $A_2$, respectively.

Summarizing, the upper bound for the nonlinear term is

(7.30) $$O(\|f\|_{l_{q,\infty}}^q \tilde{\lambda}_n^{p-q} + n^{-p\alpha/2}\log n).$$

*Rate results.* The overall rate of convergence depends on the three main contributing terms, the bias term, the linear term and the nonlinear term,

(7.31) $$\mathrm{E}\|f - \hat{f}_n\|_{L^p(g)}^p = O(\tilde{\lambda}_n^{2\delta p}) + O(n^{-p\alpha/2})$$
$$+ O(\|f\|_{l_{q,\infty}}^q \tilde{\lambda}_n^{p-q}) + O(n^{-p\alpha/2}\log n).$$

*The dense phase.* This is the region where $\alpha > \alpha_D$, $\delta = s$ and $s > (p-\pi)/2\pi$. For $\alpha > \alpha_D$, the linear term is negligible, since $n^{-p\alpha/2} = o(n^{-p\alpha_D/2})$. The bias term is negligible too since $\tilde{\lambda}_n^{2ps} = (\log n)^{2sp} n^{-sp} = o(n^{-p\alpha_D/2})$ for $s \geq 1/2$. For the nonlinear term we note that, for $q = q_D := \frac{p}{2s+1}$,

$$\tilde{\lambda}_n^{p-q} = \tilde{\lambda}_n^{2ps/(2s+1)} = \tilde{\lambda}_n^{p\alpha_D} = n^{-sp/(2s+1)}(\log n)^{2sp/(2s+1)},$$

which is the convergence rate under the dense regime. To complete the proof, we apply the Besov embedding 1 of Theorem A.1, noting that, in the dense region, we always have $\pi > q_D$.

*The sparse phase.* Here, $\alpha > \alpha_S$, $\delta = s - (1/\pi - 1/p)$ and $s < (p-\pi)/2\pi$. For $\alpha > \alpha_S$, the linear term contribution is negligible since $n^{-p\alpha/2} = o(n^{-p\alpha_S/2})$. The bias term is negligible, too, since, for $s > 1/\pi$, we have $\tilde{\lambda}_n^{2p\delta} = o(n^{-p\alpha_S/2})$. For the nonlinear term we note that for $q = q_S = \frac{p/2-1}{s+(1/2-1/\pi)}$ we have $\tilde{\lambda}_n^{p-q} = \tilde{\lambda}_n^{p\alpha_S} = n^{-p\alpha_S/2}(\log n)^{2p\alpha_S/2}$ which is the convergence rate under the sparse regime. To complete the proof, we apply the Besov embedding 3 of Theorem A.1, noting that, in the sparse region, we always have $\pi < q_D$.

*The LRD phase.* This is the region where $\alpha \leq \min(\alpha_S, \alpha_D)$. In this case, we have, for $s$ in the dense setting, $n^{-p/2\alpha_D} = o(n^{-p/2\alpha})$ and, for $s$ in the sparse setting, $n^{-p/2\alpha_S} = o(n^{-p/2\alpha})$.

7.7. *Proof of Theorem 4.8.* Let us write

$$\tilde{f}_n(x) - f(x)$$
$$= \left\{ \sum_{(j,k) \in \Lambda_1} \beta_{j,k} \psi_{j,k}(G(x)) - f(x) \right\}$$
$$+ \left\{ \sum_{(j,k) \in \Lambda_1} \tilde{\beta}_{j,k} \mathbb{I}\{|\tilde{\beta}_{j,k}| \geq \tau_0 \lambda_{n,j}\} \psi_{j,k}(G(x)) - \sum_{(j,k) \in \Lambda_1} \beta_{j,k} \psi_{j,k}(G(x)) \right\}$$
$$+ \sum_{(j,k) \in \Lambda_1} (\tilde{\beta}_{j,k} \mathbb{I}\{|\tilde{\beta}_{j,k}| \geq \tau_0 \lambda_{n,j}\} - \beta_{j,k}) \{\psi_{j,k}(\hat{G}_n(x)) - \psi_{j,k}(G(x))\}$$
$$+ \sum_{(j,k) \in \Lambda_1} \beta_{j,k} \{\psi_{j,k}(\hat{G}_n(x)) - \psi_{j,k}(G(x))\}.$$

Now, replacing Lemma 7.1 and Proposition 7.2 with Lemma 7.5, we may proceed as in the proof of Theorem 4.1 to conclude that the second part of the above decomposition is bounded with the desired rate. The third part is clearly of the smaller order than the second one. Furthermore, for the bias term we have

$$\left\| \sum_{j > j_1} \sum_{k \in D_j} \beta_{j,k} \psi_{j,k}(G(\cdot)) \right\|_{L^p(g)}^p \leq C \|f \circ G^{-1}\|_{\mathcal{B}_{p,r}^\delta}^p 2^{-j_1 \delta p}$$
$$= O((\log n / n)^{\delta p/2}).$$

Note that we have a different bound than compared to (7.28), since, here, we stopped earlier (i.e., $2^{j_1} \sim \sqrt{n/\log n}$). Nevertheless, comparing the bias term with the rate in the dense phase, we see that, with the choice $\delta = s$, we have $n^{-sp/2} < n^{-sp/(2s+1)}$, if $s > 1/2$. Furthermore, in the sparse phase, by choosing $\delta = s - (1/\pi - 1/p)$, we see that the bias is negligible as long as $s > 1/\pi + 1/2$.

Therefore, to finish the proof of Theorem 4.8, it suffices to bound the last part. We have, by using Hölder inequality,

$$\mathrm{E} \left\| \sum_{(j,k) \in \Lambda_1} \beta_{j,k} \{\psi_{j,k}(\hat{G}_n(x)) - \psi_{j,k}(G(x))\} \right\|_{L^p(g)}^p$$
$$\leq \mathrm{E} \|\hat{G}_n - G\|_\infty^p \sum_{(j,k) \in \Lambda_1} \|\beta_{j,k} \psi'_{j,k}(G(\cdot))\|_{L^p(g)}^p$$
$$= O(n^{-p/2}) \sum_{(j,k) \in \Lambda_1} 2^{j(3p/2-1)} |\beta_{j,k}|^p$$
$$= O(n^{-p/2}) \sum_{j \leq j_1} 2^{jp} 2^{-j\delta p} \left( 2^{jp(\delta+1/2-1/p)} \sum_k |\beta_{j,k}|^p \right)$$

$$= O(n^{-p/2} 2^{j_1 p(1-\delta)}) \|f \circ G^{-1}\|^p_{\mathcal{B}^\delta_{p,\infty}}$$

$$= O\left(n^{-p/2} \left(\frac{n}{\log n}\right)^{p(1-\delta)/2}\right) \|f \circ G^{-1}\|^p_{\mathcal{B}^\delta_{p,\infty}}$$

$$= O(\max\{n^{-p/2}, n^{-\delta p/2}\}) \|f \circ G^{-1}\|^p_{\mathcal{B}^\delta_{p,\infty}}.$$

If $p \leq \pi$, take $\delta = s$, so that $\mathcal{B}^s_{\pi,\infty} \subseteq \mathcal{B}^s_{p,\infty}$. The above rate is then $O(n^{-sp/(2s+1)})$ for $s \geq 1/2$. If $p > \pi$, take $\delta = s - (1/\pi - 1/p)$. The above rate is $\max\{n^{-p/2}, n^{-(s-(1/\pi-1/p))p/2}\}$ and is smaller than $n^{-p\alpha_S/2}$ as long as $s > 1/2 + 1/\pi$.

REMARK 7.6. Let us consider

$$\mathcal{F} = \{f_{j,k} = \beta_{j,k}\psi_{j,k}, j \geq 1, k = 0, \ldots, 2^j - 1\},$$

where $\beta_{j,k} = 2^{-j(s+1/2-1/\pi)}$, and we assume that $\beta_{j,k}$ are known. We recover the function $f_{j,k}$ by using the estimator $\beta_{j,k}\psi_{j,k}(\hat{G}_n(\cdot))$. Its expected weighted mean square loss, $\mathrm{E}\|\cdot\|^2_{L^2(g)}$, is

$$\beta^2_{j,k} \mathrm{E}\left[\int |\psi_{j,k}(\hat{G}_n(x)) - \psi_{j,k}(G(x))|^2 g(x)\, dx\right].$$

By considering the first term in the Taylor expansion, the above expected value is of the order

$$\mathrm{E}\left[\int \{\psi'_{j,k}(G(x))(\hat{G}_n(x) - G(x))\}^2 g(x)\, dx\right]$$

$$\sim \int \{\psi'_{j,k}(u)\}^2 u(1-u)\, du$$

$$= \frac{2^{3j}}{n} \int \{\psi'(2^j u - k)\}^2 u(1-u)\, du$$

$$= \frac{2^{2j}}{n} \int \{\psi'(v)\}^2 \left(\frac{v+k}{2^j}\right)\left(1 - \frac{v+k}{2^j}\right) dv.$$

Take $k = 2^{j/2}$. Then, the above expression is of the order $2^{3j/2}/n$. Now, if we choose $j \sim \frac{n}{\log n}$, then the expected weighted mean square loss is of the order

$$(7.32) \quad \beta^2_{j,k} 2^{j/2} \frac{2^j}{n} \sim 2^{-2j(s+1/2-1/\pi)} 2^{j/2} \frac{1}{\log n} \sim n^{-2(s+1/4-1/\pi)} \times \log \text{ term}.$$

Choose, for simplicity, $\pi = 1$. Since also $p = 2$, there is no sparse phase and the only restriction (in the Theorem 4.1) in the dense phase is $s > 1$. However, we note that the rate in (7.32) is of the smaller order than $n^{-2s/(2s+1)}$ if and only if $s < \frac{3}{8} - \frac{1}{8}\sqrt{33}$ or $s \geq \frac{3}{8} + \frac{1}{8}\sqrt{33} > 1$. Consequently, we cannot stop the fully adaptive estimator at $2^{j_1} \sim \frac{n}{\log n}$ and keep the same restriction on $s$, as in the case of the partially adaptive one.

## APPENDIX: BESOV EMBEDDING IN $L_{Q,\infty}$ SPACES

We give a simplified version of Theorem 6.2 [18], when the dimension $d = 1$ and $\sigma_j = 1$. Let $\mu$ will denote the measure such that for $j \in \mathbb{N}, k \in \mathbb{N}$,

$$\text{(A.1)} \quad \mu\{(j,k)\} = \|\psi_{j,k}\|_p^p = 2^{j(p/2-1)}\|\psi\|_p^p,$$

$$\text{(A.2)} \quad l_{q,\infty} := \left\{ f = \sum_{j,k} \beta_{j,k}\psi_{j,k}, \right.$$

$$\left. \|f\|_{l_{q,\infty}} := \sup_{\lambda > 0} \lambda^q \mu\{(j,k) : |\beta_{j,k}| > \lambda\} < \infty \right\}$$

and

$$l_q := \left\{ f = \sum_{j,k} \beta_{j,k}\psi_{j,k} \in L^p, \|f\|_{l_q} := \left( \sum_{j,k \in A_j} |\beta_{j,k}|^q \mu\{(j,k)\} \right)^{1/q} < \infty \right\},$$

where $A_j$ is a set of cardinality proportional to $2^j$.

THEOREM A.1. *Let $0 < p < \infty, 0 \leq s \leq \infty$ be fixed and let $q_D = p/(2s+1)$:*

1. *If $\pi > q_D$, then for all $r$, $0 < r \leq \infty$, $\mathcal{B}_{\pi,r}^s \subset \mathcal{B}_{\pi,\infty}^s \subset l_{q_D,\infty}$.*
2. *If $\pi = q_D$, then for all $r$, $0 < r \leq \pi$, $\mathcal{B}_{\pi,r}^s \subset \mathcal{B}_{\pi,\pi}^s \subset l_\pi$. Moreover for $r > \pi$, we have:*

   - *If $p = 2$ then $\mathcal{B}_{\pi,r}^s \subset l_\pi$.*
   - *If $p > 2$ then for all $r > p$, $\mathcal{B}_{\pi,r}^s \subset \mathcal{B}_{\pi,\infty}^s \subset l_r$.*

3. *If $2/(2s+1) < \pi < q_D$, for all $0 < r \leq \infty$, $\mathcal{B}_{\pi,r}^s \subset \mathcal{B}_{\pi,\infty}^s \subset l_{q_S,\infty}$, where $q_S = \frac{p/2-1}{s+(1/2-1/\pi)}$.*

DEPARTMENT OF MATHEMATICS
AND STATISTICS
UNIVERSITY OF OTTAWA
585 KING EDWARD AVENUE
OTTAWA ON K1N 6N5
CANADA
E-MAIL: rkulik@uottawa.ca